# MAXIMA OF ASYMPTOTICALLY GAUSSIAN RANDOM FIELDS AND MODERATE DEVIATION APPROXIMATIONS TO BOUNDARY CROSSING PROBABILITIES OF SUMS OF RANDOM VARIABLES WITH MULTIDIMENSIONAL INDICES

By Hock Peng Chan[1] and Tze Leung Lai[2]

*National University of Singapore and Stanford University*

Several classical results on boundary crossing probabilities of Brownian motion and random walks are extended to asymptotically Gaussian random fields, which include sums of i.i.d. random variables with multidimensional indices, multivariate empirical processes, and scan statistics in change-point and signal detection as special cases. Some key ingredients in these extensions are moderate deviation approximations to marginal tail probabilities and weak convergence of the conditional distributions of certain "clumps" around high-level crossings. We also discuss how these results are related to the Poisson clumping heuristic and tube formulas of Gaussian random fields, and describe their applications to laws of the iterated logarithm in the form of the Kolmogorov–Erdős–Feller integral tests.

**1. Introduction.** The goal of this paper is to extend a number of classical results on boundary crossing probabilities of Brownian motion and random walks to much more general stochastic processes involving multidimensional indices (i.e., random fields). These extensions were motivated by applications to signal detection and change-point problems; see Example 2.2 and the last two paragraphs of Section 4. Other applications include the laws of the iterated logarithm for sums of i.i.d. random variables with multidimensional indices (see Section 3), Kolmogorov–Smirnov statistics of

Received April 2004; revised November 2004.
[1]Supported by the National University of Singapore.
[2]Supported by the National Science Foundation and the Institute of Mathematical Research at University of Hong Kong.

multivariate distributions and sums of linear processes with long-range dependence (see Section 4). To begin with, let $\{W(t): t \geq 0\}$ be Brownian motion and let $T_c = \inf\{t \geq 0 : W(t) \geq b_c(t)\}$ be the first time when Brownian motion crosses a positive continuously differentiable boundary $b_c$. Strassen [34], Jennen and Lerche [22], Wichura [37] and others have shown that $T_c$ has a density function $p_c$ and that under certain additional conditions on $b_c$, $p_c$ has the "tangent approximation"

$$(1.1) \qquad p_c(t) \doteq t^{-3/2} a_c(t) \varphi(b_c(t)/\sqrt{t}),$$

where $\varphi(x) = (2\pi)^{-1/2} e^{-x^2/2}$ is the standard normal density function and $a_c(t) = b_c(t) - tb'_c(t)$. Note that in the case of a linear boundary $b_c(t) = a + \beta t$ (with $a > 0$ and $\beta > 0$), the well-known Bachelier–Lévy formula yields $p_c(t) = t^{-3/2} a \varphi(b_c(t)/\sqrt{t})$, so (1.1) simply replaces $\alpha$ by the intercept $\alpha_c(t)$ of the tangent line passing through $(t, b_c(t))$, and is therefore called a "tangent approximation." For concave boundaries $b_c(t) = b(t)$ that become infinite as $t \to \infty$, one typically has $b'(t) = o(b(t)/t)$, so one can replace $a_c(t)$ in (1.1) by $b(t)$. There is a close connection between this approximation to $p_c(t)$ and the Kolmogorov–Erdős–Feller test, which yields for nondecreasing $b(t)/\sqrt{t}$ the 0–1 dichotomy

$(1.2) P\{W(t) < b(t) \text{ for all large } t\} = 1 \text{ (or 0)} \qquad \text{if } \mathcal{I}(b) < \infty \text{ (or } = \infty),$

where $\mathcal{I}(b) = \int_1^\infty t^{-3/2} b(t) \varphi(b(t)/\sqrt{t}) \, dt < \infty$. Similarly, if $S_n = X_1 + \cdots + X_n$ with $EX_1 = 0$, $EX_1^2 = 1$ and $E|X_1|^3 < \infty$, then for all $n \geq 1$,

$(1.3) \ P\{S_n < b(n) \text{ for all large } n\} = 1 \text{ (or 0)} \qquad \text{if } \mathcal{I}(b) < \infty \text{ (or } = \infty).$

If we think of the random walk $\{S_n, n \geq n_c\}$ in (1.3) as an "asymptotic" Brownian motion as $n_c \to \infty$, then (1.3) can be regarded as the generalization of (1.2) to processes that behave like Brownian motion. This suggests that if (1.1) and (1.2) can be extended to more general Gaussian processes, then they may even be expected to hold much more generally for processes that are "asymptotically Gaussian." In view of the functional central limit theorem for sums of weakly dependent or long-memory random variables, the scope of applications of such results would be very broad. Unfortunately, functional central limit theorems, which are about the "central" part of the limiting Gaussian distributions, are not the right tools to handle the "rare" events in the high-level crossings as in (1.1) and (1.3).

To extend (1.1) and (1.2) to much more general processes, our approach uses (i) moderate deviation approximations to marginal tail probabilities and (ii) weak convergence (to a limiting Gaussian process) of a certain conditional process given that the process attains a high level near the boundary at time $t$. Another key idea of our extension is to relax the requirement that the left-hand side of (1.1) be a first exit density. Instead we regard it as a "local" exit density at time $t$ so that the probability that the process ever crosses the boundary within time interval $D$ is asymptotically equal to the integral of the right-hand side of (1.1) over $D$. Not only does this avoid the technical assumptions that need to be imposed to ensure that the first exit time $T_c$ indeed has a density with respect to Lebesgue measure, but it also dispenses with the notion of having a well-ordered set $D$ so that the "first" time of exit can be defined. This enables us to extend our approach to random fields (with multidimensional time that is not well ordered). Section 2 gives basic assumptions for these "asymptotically Gaussian" random fields and states the main theorems that provide generalizations of (1.1) and (1.3). Applying these theorems to Gaussian random fields yields new results in Theorem 2.1 for the maxima of Gaussian random fields. Section 5 gives the proofs. Connections to Aldous' [4] Poisson clumping heuristic and the Hotelling–Weyl tube formulas are also discussed in Section 2.

**2. Basic results and discussion.** We begin with some notation that will be used throughout the paper. Let $\psi(c) = (2\pi c^2)^{-1/2} \exp(-c^2/2)$. For vectors $t, u \in \mathbf{R}^d$, the relation $t \leq u$ means $t_i \leq u_i$ for all $i$ and $t < u$ means $t_i < u_i$ for all $i$. Also $\lfloor \cdot \rfloor$ will be used to denote the greatest integer function, $\|\cdot\|$ the (Euclidean) norm of a vector, $|\cdot|$ the determinant of a square matrix and $v(\cdot)$ the $d$-dimensional volume (or content) of a Jordan measurable set. For $\zeta > 0$, let

$$I_{t,\zeta} = \prod_{i=1}^{d} [t_i, t_i + \zeta).$$

For $D \subset \mathbf{R}^d$ and $\delta > 0$, define $[D]_\delta = \{t + u : t \in D, \|u\| < \delta\}$. We shall also use $\nabla$ and $\nabla^2$ to denote the gradient vector and Hessian matrix, respectively, of a function. Let $\mathbf{S}^{d-1}$ denote the $(d-1)$-dimensional unit sphere, and let $\mathbf{Z}_+$ ($\mathbf{R}_+$) denote the set of positive integers (real numbers). Let $0 < \alpha \leq 2$ and let $\{W_t(u) : u \in [0,\infty)^d\}$ be a continuous Gaussian random field (whose continuity follows from Theorem 2.1 of [25]) such that

$$W_t(0) = 0,$$
$$E[W_t(u)] = -\|u\|^\alpha r_t(u/\|u\|)/2,$$
(2.1)
$$\mathrm{Cov}(W_t(u), W_t(v)) = [\|u\|^\alpha r_t(u/\|u\|) + \|v\|^\alpha r_t(v/\|v\|) - \|u-v\|^\alpha r_t((u-v)/\|u-v\|)]/2,$$

where $r_t : \mathbf{S}^{d-1} \to \mathbf{R}_+$ is a continuous function satisfying

(2.2) $$\sup_{v \in \mathbf{S}^{d-1}} |r_t(v) - r_u(v)| \to 0 \quad \text{as } u \to t.$$

Of particular importance in the subsequent development are

$$H_K(t) = \int_0^\infty e^y P\left\{\sup_{0 \leq u_i \leq K \ \forall i} W_t(u) > y\right\} dy,$$

(2.3)

$$H(t) = \lim_{K \to \infty} K^{-d} H_K(t),$$

which are shown to be well defined in Theorems 2.4 and 2.5.

Let $X$ be a stationary, isotropic Gaussian random field such that $EX(0) = 0$, $EX^2(0) = 1$ and

(2.4) $\quad E[X(0)X(u)] = 1 - (1 + o(1))\|u\|^\alpha L(\|u\|) \quad$ as $u \to 0$,

for some $0 < \alpha \leq 2$ and slowly varying function $L$. Let

(2.5) $\quad \Delta_c = \min\{x > 0 : x^\alpha L(x) = (2c^2)^{-1}\}.$

For example, if $L(x) \equiv 1$, then $\Delta_c = (2c^2)^{-1/\alpha}$. Let $D$ be a bounded, Jordan measurable set such that $[D]_\delta$ lies in the domain of $X$ for some $\delta > 0$. Then by Theorem 2.1 of [31],

(2.6) $\quad P\left\{\sup_{t \in D} X(t) > c\right\} \sim \psi(c) \Delta_c^{-d} v(D) H,$

where $H = \lim_{K \to \infty} K^{-d} \int_0^\infty e^y P\{\sup_{0 \leq u_i \leq K \ \forall i} W_0(u) > y\} dy$ is a positive, finite constant and $W_0$ is the Gaussian random field defined in (2.1) with $r_0(u) \equiv 1$. Our goal is to extend (2.6) first to more general Gaussian random fields satisfying

(2.7) $\ E[X(t)X(t+u)] = 1 - (1 + o(1))\|u\|^\alpha L(\|u\|) r_t(u/\|u\|) \quad$ as $u \to 0$,

uniformly over $t \in [D]_\delta$. We then extend (2.6) to non-Gaussian random fields that are asymptotically Gaussian in a moderate deviation sense.

2.1. *Gaussian random fields.* Let $X$ be a Gaussian random field such that $EX(t) = 0$, $EX^2(t) = 1$ for all $t$. Let $D$ be such that $[D]_\delta$ is a subset of the domain of $X$ for some $\delta > 0$. The following theorem, whose proof is given in Section 5, generalizes (2.6) far beyond the stationary isotropic framework considered by Qualls and Watanabe [31] under (2.4).

THEOREM 2.1. *Suppose the Gaussian random field $X$ satisfies condition (2.7), in which $0 < \alpha \leq 2$ and $r_t : \mathbf{S}^{d-1} \to \mathbf{R}_+$ is a continuous function such that the convergence in (2.2) is uniform in $t \in [D]_\delta$ and $\sup_{t \in [D]_\delta, v \in \mathbf{S}^{d-1}} r_t(v) < \infty$. Then, with $H(t)$ defined by (2.3),*

(2.8) $\quad P\left\{\sup_{u \in I_{t,\ell_c \Delta_c}} X(u) > c\right\} \sim \ell_c^d \psi(c) H(t)$

uniformly over $t \in D$, as $c \to \infty$ and $\ell_c \to \infty$ with $\ell_c = o(\Delta_c^{-1})$. Moreover, if $D$ is bounded and Jordan measurable, then as $c \to \infty$,

$$(2.9) \qquad P\left\{\sup_{t \in D} X(t) > c\right\} \sim \psi(c) \Delta_c^{-d} \int_D H(t) \, dt.$$

The following special case of Theorem 2.1, with $d = 2$, demonstrates the usefulness of including the function $r_t$ on $\mathbf{S}^{d-1}$ in (2.2) when (2.1) is extended to nonstationary Gaussian random fields. It will be discussed further in Example 2.10 and at the end of Section 4.

EXAMPLE 2.2. Let $X(t_1, t_2) = (t_2 - t_1)^{-1/2}[W(t_2) - W(t_1)]$, where $W(\cdot)$ is Brownian motion, and $D = \{(t_1, t_2) : 0 \le t_1 < t_2 \le a, a_1 \le t_2 - t_1 \le a_2\}$ with $0 < a_1 < a_2 < a$. Then

$$E[X(t)X(t+u)] = 1 - (1 + o(1)) \frac{|u_1| + |u_2|}{2(t_2 - t_1)},$$

as $u \to 0$. Hence (2.7) is satisfied with $\alpha = 1$, $L(\|u\|) \equiv 1$ and $r_t(u) = (|u_1| + |u_2|)/[2(t_2 - t_1)]$. Therefore $\Delta_c = (2c^2)^{-1}$ in view of (2.5), and $H(t) = 2^{-4}(t_2 - t_1)^{-2}$ by Lemma 2.3 below. Application of Theorem 2.1 then yields that as $c \to \infty$,

$$P\{(t_2 - t_1)^{-1/2}[W(t_2) - W(t_1)] > c$$

$$\text{for some } 0 \le t_1 < t_2 \le a, a_1 \le t_2 - t_1 \le a_2\}$$

$$(2.10) \qquad \sim \psi(c)(2c^2)^2 \int_D 2^{-4}(t_2 - t_1)^{-2} \, dt_1 \, dt_2$$

$$= \psi(c)(c^4/4) \int_{a_1}^{a_2} \int_0^{a-s} s^{-2} \, dt_1 \, ds$$

$$= \psi(c)(c^4/4)[a(a_1^{-1} - a_2^{-1}) - \log(a_2/a_1)].$$

LEMMA 2.3. Let $\{W_t(u) : u \in [0, \infty)^d\}$ be a continuous Gaussian random field such that for some positive functions $\beta_1, \ldots, \beta_d$,

$$E[W_t(u)] = -\sum_{i=1}^d \beta_i(t) u_i / 2$$

and

$$\operatorname{Cov}[W_t(u), W_t(v)] = \sum_{i=1}^d \beta_i(t)(u_i + v_i - |u_i - v_i|)/2 = \sum_{i=1}^d \beta_i(t) \min(u_i, v_i).$$

Then $H(t) = 2^{-d} \prod_{i=1}^d \beta_i(t)$.

PROOF. For $u \geq 0$, $W_t(u) = \sum_{i=1}^d B_{i,t}(u_i)$, where $\{B_{i,t}\}_{1 \leq i \leq d}$ are independent Gaussian processes with independent increments, $E[B_{i,t}(u_i)] = -\beta_i(t)u_i/2$ and $\text{Var}(B_{i,t}(u_i)) = \beta_i(t)u_i$, so $B_{i,t}(u_i) \stackrel{\mathcal{L}}{=} W(\beta_i(t)u_i) - \beta_i(t)u_i/2$. As $K \to \infty$,

$$H_K(t) = \int_0^\infty e^y P\left\{\sum_{i=1}^d \sup_{0 \leq u_i \leq K} B_{i,t}(u_i) > y\right\} dy$$

$$= \int_0^\infty (e^y - 1) P\left\{\sum_{i=1}^d \sup_{0 \leq u_i \leq K} B_{i,t}(u_i) \in dy\right\}$$

$$\sim E\left\{\exp\left[\sum_{i=1}^d \sup_{0 \leq u_i \leq K} B_{i,t}(u_i)\right]\right\}$$

$$= \prod_{i=1}^d E\left\{\exp\left(\sup_{0 \leq u_i \leq K} [W(\beta_i(t)u_i) - \beta_i(t)u_i/2]\right)\right\},$$

so $H_K(t) \sim \prod_{i=1}^d [\beta_i(t)K/2]$; see [18], (1.8.11) for the last asymptotic relation. □

2.2. *Asymptotically Gaussian random fields.* Theorem 2.1 is derived in Section 5 as a special case of a more general result on *asymptotically Gaussian* random fields satisfying conditions (C) and (A1)–(A5) below. Specifically, for $c > 0$, let $X_c$ be random fields such that $EX_c(t) = 0$, $EX_c^2(t) = 1$ for all $c$ and $t$. Let $D$ be such that $[D]_\delta$ is a subset of the domain of $X_c$ for some $\delta > 0$ and all $c$ large enough. Define $\rho_c(t, u) = E[X_c(t)X_c(u)]$. In analogy with (2.7), assume that there exist $0 < \alpha \leq 2$ and a slowly varying function $L$ such that as $u \to 0$,

(C)  $\qquad \rho_c(t, t+u) = 1 - (1+o(1))\|u\|^\alpha L(\|u\|) r_t(u/\|u\|)$

uniformly over $t \in [D]_\delta$ and compact sets of $u/\Delta_c > 0$. Moreover, assume that the following conditions also hold uniformly over $t \in [D]_\delta$, as $c \to \infty$:

(A1)  $\qquad P\{X_c(t) > c - y/c\} \sim \psi(c - y/c)$

uniformly over positive, bounded values of $y$. The convergence in (2.2) is assumed to be uniform in $t \in [D]_\delta$, with $\sup_{t \in [D]_\delta, v \in \mathbf{S}^{d-1}} r_t(v) < \infty$. Moreover, for any $a > 0$ and positive integers $m$, as $c \to \infty$,

(A2)  $\qquad \begin{aligned} &\{c[X_c(t + ak\Delta_c) - X_c(t)] : 0 \leq k_i < m\} | X_c(t) = c - y/c \\ &\Rightarrow \{W_t(ak) : 0 \leq k_i < m\} \end{aligned}$

uniformly over positive, bounded values of $y$, where we use "$|X_c(t) = c - y/c$" to denote that the distribution is conditional on $X_c(t) = c - y/c$. In addition, there exists a positive function $h$ such that $\lim_{y \to \infty} h(y) = 0$ and

(A3) $\quad P\{X_c(t + u\Delta_c) > c - \gamma/c, X_c(t) \leq c - y/c\} \leq h(y)\psi(c)$

for all $u \geq 0$ and $\gamma > 0$, and there exist nonincreasing functions $N_a$ on $\mathbf{R}_+$ and positive constants $\gamma_a$ such that $\gamma_a \to 0$ and $N_a(\gamma_a) + \int_1^\infty \omega^s N_a(\gamma_a + \omega) \, d\omega = o(a^d)$ as $a \to 0$, and

(A4) $\quad P\left\{\sup_{0 \leq u \leq a} X_c(t + u\Delta_c) > c, X_c(t) \leq c - \gamma/c\right\} \leq N_a(\gamma)\psi(c),$

for all $\gamma_a \leq \gamma \leq c$ and $s > 0$.

Whereas (A1) refers to the marginal distribution of $X_c(t)$, saying that $\{X_c(t) > c - y/c\}$ has probability like that of a standard normal, the joint distribution of $X_c(\cdot)$ is assumed in (A2) to be asymptotically normal in the sense of weak convergence for local increments conditioned on $X_c(t) = c - y/c$. Note that the same $\alpha, L(\cdot)$ and $r_t(\cdot)$ appear in (C) and the mean and covariance functions (2.1) of the Gaussian field $W_t(\cdot)$ in (A2). In fact, if $X_c = X$ is a Gaussian field satisfying condition (C), then (A2) holds; see the proof of Theorem 2.1 in Section 5. Assumptions (A3) and (A4) are mild technical conditions under which the probability of $\sup_{u \in I_{t, K\Delta_c}} X_c(u)$ exceeding $c$ can be computed via (A1) and (A2) after the cube $I_{t, K\Delta_c} = \prod_{i=1}^d [t_i, t_i + K\Delta_c)$ is discretized by the grid points $t + ka\Delta_c$ ($0 \leq k_i < m$) with $a = K/m$, leading to the following.

THEOREM 2.4. *Let $K > 0$. Assume* (C) *and* (A1)–(A4). *Then as $c \to \infty$,*

$$P\left\{\sup_{u \in I_{t, K\Delta_c}} X_c(u) > c\right\} \sim \psi(c)[1 + H_K(t)]$$

*uniformly over $t \in [D]_\delta$, where $H_K(t)$ is defined in* (2.3) *and is finite and uniformly continuous in $t \in [D]_\delta$.*

To derive an analogue of (2.6) for $P\{\sup_{t \in D} X_c(t) > c\}$ in which $X_c$ satisfies (C) and (A1)–(A4), we can sum the asymptotic formula in Theorem 2.4 over $t \in (K\Delta_c \mathbf{Z})^d \cap D$ if the joint occurence of two events (associated with two such cubes) is negligible in comparison with the probability associated with a single cube. The following simple condition ensures this: There exists a nonincreasing function $f: [0, \infty) \to \mathbf{R}_+$ such that $f(\|r\|) = O(e^{-\|r\|^p})$ for some $p > 0$ and for all $\gamma > 0$ and $c$ sufficiently large,

(A5) $\quad P\{X_c(t) > c - \gamma/c, X_c(t + u\Delta_c) > c - \gamma/c\} \leq \psi(c - \gamma/c)f(\|u\|)$

uniformly in $t$ and $t + u\Delta_c$ belonging to $[D]_\delta$.

THEOREM 2.5. *Assume* (C) *and* (A1)–(A5). *Then as* $c \to \infty$ *and* $\ell_c \to \infty$ *such that* $\ell_c = o(\Delta_c^{-1})$,

$$P\left\{\sup_{u \in I_{t,\ell_c\Delta_c}} X_c(u) > c\right\} \sim \ell_c^d \psi(c) H(t), \tag{2.11}$$

$$P\left\{\sup_{u \in I_{t,\ell_c\Delta_c}} X_c(u) > c, \sup_{v \in B \setminus I_{t,\ell_c\Delta_c}} X_c(v) > c\right\} = o(\ell_c^d \psi(c)), \tag{2.12}$$

*uniformly over* $t \in D$ *and over subsets* $B$ *of* $[D]_\delta$ *with bounded volume, where* $H(t)$ *is defined in* (2.3) *and is uniformly continuous and bounded below on* $D$.

Dividing (2.11) by $(\ell_c \Delta_c)^d$, which is the volume of $I_{t,\ell_c\Delta_c}$, yields an asymptotic boundary crossing "density" $\Delta_c^{-d} \psi(c) H(t)$ of $X_c$ at $t$. By integrating this "density" over $D$, or more precisely, by summing (2.11) over the "tiles" $I_{t,\ell_c\Delta_c}$ of $D$ and applying (2.12) together with the fact that $D$ is bounded and Jordan measurable, we obtain the following generalization of the Qualls–Watanabe result (2.6) on stationary isotropic Gaussian random fields.

COROLLARY 2.6. *Assume* (C) *and* (A1)–(A5). *Let* $D$ *be a bounded, Jordan measurable set. Then*

$$P\left\{\sup_{t \in D} X_c(t) > c\right\} \sim \psi(c) \Delta_c^{-d} \int_D H(t)\,dt \qquad \text{as } c \to \infty. \tag{2.13}$$

We can extend Corollary 2.6 to sets $D_c$ that grow with $c$. The assumption that $D$ be Jordan measurable [i.e., for any $\varepsilon > 0$, the boundary $\partial D$ of $D$ can be covered by rectangles $U_1, U_2, \ldots$ such that $\sum_{i=1}^\infty v(U_i) < \varepsilon$] and bounded in Corollary 2.6 is used to show that $\sum_{t \in (\zeta \mathbf{Z})^d, I_{t,\zeta} \cap \partial D \neq \varnothing} v(I_{t,\zeta}) \to 0$ as $\zeta \to 0$. When working with sets $D_c$ that need not be bounded, we need to impose a more direct assumption (2.14) on the contribution of $\partial D_c$ to the Riemann sum. Moreover, by condition (C) or (A1)–(A5), we now mean that it holds uniformly over $t$ belonging to $[D_c]_\delta$.

COROLLARY 2.7. *Assume* (C), (A1)–(A5) *and that*

$$\sup_{t,u \in D_c} \|u - t\| = O(c^\kappa) \quad \text{and} \quad \zeta_c^d \sum_{t \in (\zeta_c \mathbf{Z})^d, I_{t,\zeta_c} \cap \partial D_c \neq \varnothing} H(t) = o(v(D_c)) \tag{2.14}$$

*for some* $\kappa > 0$ *and positive* $\zeta_c$ *with* $\zeta_c \to 0$ *and* $\Delta_c = o(\zeta_c)$. *Then as* $c \to \infty$,

$$P\left\{\sup_{t \in D_c} X_c(t) > c\right\} \sim \psi(c) \Delta_c^{-d} \int_{D_c} H(t)\,dt. \tag{2.15}$$

2.3. *Boundary crossing probabilities.* To extend the conclusion of Corollary 2.7 to the boundary crossing probability $P\{X_c(t) > b_c(t)$ for some $t \in D_c\}$, we proceed similarly by using the probabilities $p_c(t) = P\{X_c(s) > b_c(s)$ for some $s \in I_{t,\zeta_c}\}$ as building blocks, where $\zeta_c \to 0$ is so chosen that

$$(2.16) \quad \sup_{t \in [D_c]_\delta} \Delta_{b_c(t)} = o(\zeta_c) \quad \left(\text{hence} \inf_{t \in [D_c]_\delta} b_c(t) \to \infty\right) \quad \text{as } c \to \infty.$$

Whereas (A1)–(A5) are related to the time-invariant boundary $c$ to be crossed by $X_c(\cdot)$, we can formulate similar assumptions when $c$ is replaced by a time-varying boundary $b_c(\cdot)$. Let

$$(2.17) \quad \underline{b}_c = \inf_{u \in [D_c]_\delta} b_c(u), \qquad \overline{b}_c = \sup_{u \in [D_c]_\delta} b_c(u).$$

Analogous to (A1)–(A5), assume that the following conditions hold, as $c \to \infty$, uniformly in $t \in [D_c]_\delta$ and $\underline{b}_c/2 \leq z \leq \overline{b}_c$:

(B1) $$P\{X_c(t) > z\} \sim \psi(z),$$

(B2) $$\{z[X_c(t + ak\Delta_z) - X_c(t)] : 0 \leq k_i < m\} | X_c(t) = z - \gamma/z$$
$$\Rightarrow \{W_t(ak) : 0 \leq k_i < m\},$$

for any $a > 0$ and positive integers $m$, the convergence being uniform over positive, bounded values of $\gamma$; moreover, the convergence in (2.2) is assumed to be uniform in $t \in [D_c]_\delta$, with $\sup_{t \in [D_c]_\delta, v \in \mathbf{S}^{d-1}} r_t(v) < \infty$. In addition, there exists a positive function $h$ such that $\lim_{y \to \infty} h(y) = 0$ and

(B3) $$P\{X_c(t + u\Delta_z) > z - \gamma/z, X_c(t) \leq z - y/z\} \leq h(y)\psi(z)$$

for all $u \geq 0$ and $\gamma > 0$, and there exist nonincreasing functions $N_a$ on $\mathbf{R}_+$ and positive constants $\gamma_a$ such that $\gamma_a \to 0$ and $N_a(\gamma_a) + \int_1^\infty \omega^s N_a(\gamma_a + \omega) \, d\omega = o(a^d)$ as $a \to 0$, and

(B4) $$P\left\{\sup_{0 \leq u \leq a} X_c(t + u\Delta_z) > z, X_c(t) \leq z - \gamma/z\right\} \leq N_a(\gamma)\psi(z),$$

for all $\gamma_a \leq \gamma \leq z$ and $s > 0$. Moreover, there exists a nonincreasing function $f:[0,\infty) \to \mathbf{R}_+$, with $f(\|r\|) = O(e^{-\|r\|^p})$ for some $p > 0$, such that for $\gamma > 0$ and $c$ sufficiently large

(B5) $P\{X_c(t) > z - \gamma/z, X_c(t + u\Delta_z) > z - \gamma/z\} \leq \psi(z - \gamma/z)f(\|u\|)$

uniformly in $t$ and $t + u\Delta_z$ belonging to $[D_c]_\delta$.

THEOREM 2.8. *Assume* (C) *and* (B1)–(B5). *Suppose that* (2.14) *and* (2.16) *hold for some* $\kappa > 0$ *and* $\zeta_c \to 0$ *and that*

$$(2.18) \quad \sup_{t \in [D_c]_{2\zeta_c}} [\overline{b}_c^2(t) - \underline{b}_c^2(t)] = o(1)$$

*where* $\overline{b}_c(t) = \sup_{u \in I_{t,\zeta_c}} b_c(u), \underline{b}_c(t) = \inf_{u \in I_{t,\zeta_c}} b_c(u).$

Then $P\{X_c(t) > b_c(t) \text{ for some } t \in D_c\} \sim \int_{D_c} \psi(b_c(t))\Delta_{b_c(t)}^{-d} H(t)\,dt$ as $c \to \infty$.

The next corollary specializes Theorem 2.8 to the case in which $b_c(t) = cb(t)$ for some positive function $b$ possessing continuous second derivatives on $[D]_\delta$, where $D$ is a compact Jordan measurable set. Let $b_D = \inf_{t \in D} b(t)$ and assume that $\mathcal{M} = \{t \in D : b(t) = b_D\}$ is a $q$-dimensional manifold (with boundary) such that $v_q(\mathcal{M} \cap \partial D) = 0$, in which $v_q$ denotes the $q$-dimensional volume element of the manifold. Let $T\mathcal{M}^\perp(t)$ denote the normal space of the manifold $\mathcal{M}$ at $t$. Letting $\{e_1(t), \ldots, e_{d-q}(t)\}$ be an orthonormal basis of $T\mathcal{M}^\perp(t)$, define the $d \times (d-q)$ matrix $A(t) = (e_1(t) \cdots e_{d-q}(t))$ and assume that $\nabla_\perp^2 b(t) := A'(t)\nabla^2 b(t) A(t)$ is a positive definite $q \times q$ matrix for all $t \in \mathcal{M}$.

COROLLARY 2.9. *Suppose* (C) *and* (B1)–(B5) *are satisfied with* $\alpha < 2$ *and* $D_c = D$, *a compact Jordan measurable set. Then as* $c \to \infty$,

$$P\{X_c(t) > cb(t) \text{ for some } t \in D\} \quad (2.19)$$
$$\sim \psi(cb_D) b_D^{2d/\alpha} \Delta_c^{-d} (2\pi/c^2 b_D)^{(d-q)/2} \int_\mathcal{M} |\nabla_\perp^2 b(t)|^{-1/2} H(t) v_q(dt).$$

EXAMPLE 2.10. Let $X(t_1, t_2) = (t_2 - t_1)^{-1/2}[W(t_2) - W(t_1)]$ and $X_c = X$ as in Example 2.2, where $W(\cdot)$ is Brownian motion, and let $b_c(t_1, t_2) = [c^2 + 2\log(t_2 - t_1)^{-\beta}]^{1/2}$ for some $\beta > 1$. Let $D = \{(t_1, t_2) : 0 \le t_1 < t_2 \le a, a_1 \le t_2 - t_1 \le a_2\}$ where $0 < a_1 < a_2 \le a$. Arguments similar to those used to prove Theorem 2.1 in Section 5 can be used to show that (B1)–(B5) hold uniformly in $t \in [D]_\delta$ and $\underline{b}_c/2 \le z \le \overline{b}_c$. Therefore by Lemma 2.3 and Theorem 2.8,

$$P\{(t_2 - t_1)^{-1/2}[W(t_2) - W(t_1)] > [c^2 + 2\log(t_2 - t_1)^{-\beta}]^{1/2}$$
$$\text{for some } 0 \le t_1 < t_2 \le a, a_1 \le t_2 - t_1 \le a_2\}$$
$$\sim \psi(c)(2c^2)^2 \int_D 2^{-4} e^{\beta \log(t_2 - t_1)} (t_2 - t_1)^{-2}\,dt_1\,dt_2$$
$$= \frac{c^4 \psi(c)}{4} \int_{a_1}^{a_2} \int_0^{a-s} s^{-2+\beta}\,dt_1\,ds$$
$$= \frac{c^4 \psi(c)}{4} \left( \frac{a(a_2^{\beta-1} - a_1^{\beta-1})}{\beta - 1} - \frac{(a_2^\beta - a_1^\beta)}{\beta} \right).$$

2.4. *Discussion and related literature.* Our formulation of "asymptotically Gaussian" random fields bears some resemblance to Aldous's [4] Poisson clumping heuristic, which involves i.i.d. clumps of high-level excursions of a stochastic process $X(t)$, with the stochastic structure of the clump determined by the conditional limiting process [like that in (A2)] of normalized local increments. Whereas the Poisson clumping heuristic only suggests an asymptotic approximation $P\{\sup_{t\in D} X_c(t) \leq c\}$ of the form $e^{-p_c}$ with $p_c \to 0$, our approach actually gives a rigorous derivation of an asymptotic formula for $p_c$. Instead of a single stochastic process $X(t)$, our formulation involves a family of random fields $X_c(t)$ with $EX_c(t) = 0$ and $\text{Var}(X_c(t)) = 1$. It consists of two basic components: (i) a normal approximation to the probability of $X_c(t)$ exceeding some high level (depending on $c$) in (A1) or (B1), and (ii) the weak convergence of the finite-dimensional distributions of the local increments conditioned on $X_c(t) = c - y/c$ in (A2) [or (B2)]. The covariance structure of the local increments given by condition (C) and the closely related mean and covariance functions (2.1) of the limiting Gaussian random field in (A2) [or (B2)] provide the key ingredients in the asymptotic formulas in Corollaries 2.6, 2.7 and Theorem 2.8. Theorem 2.1 and its proof show that these asymptotic formulas are the same as in the special case $X_c = X$, a zero-mean Gaussian random field satisfying condition (C). These asymptotic formulas are derived by adding up corresponding results for small cubes in (2.11), making use of (2.12) to justify the additivity.

Conditions of the type (A2) were introduced by Berman ([6], Theorem 5.1) for asymptotic approximations (as $c \to \infty$) to the probability $P\{\sup_{0 \leq t \leq T} X(t) > c\}$ of a stationary process $X(t)$ (with $d = 1$) such that $X(0)$ belongs to the domain of attraction of an extreme value distribution; see [6], Theorem 14.1 and [3], Theorem 1. We consider here general $d$, extend $X(t)$ to $X_c(t)$ and remove the stationary assumption, but restrict the limiting distribution in (A2) to be Gaussian and the marginal probabilities $P\{X_c(t) > c - y/c\}$ to be asymptotically normal. It will be shown in Sections 3 and 4 that extending a single stationary process $X(t)$ to a family of possibly nonstationary random fields $X_c(t)$ and generalizing the threshold $c$ to a moving boundary $b_c(t)$ greatly broaden the scope of applications. Some of the difficulties in proving these extensions to the nonstationary setting are explained in Remark 5.1.

Corollary 2.9 and its proof in Section 5 reveal similarities and differences between our approach and the tube formulas of Hotelling [21] and Weyl [36] whose applications to the maxima of Gaussian random fields are reviewed in Section 6 of [1]. As in [8], the use of the tubular neighborhood $U_{\xi_c}$ of the extremal manifold $\mathcal{M}$ in the proof of Corollary 2.9 is related to Laplace's method for asymptotic evaluation of the integral $\int_{D_c} \psi(b_c(t)) \Delta_{b_c(t)}^{-d} H(t) \, dt$, in which the integrand can be regarded as an "asymptotic density" of crossing the boundary $b_c$ by $X_c$ at $t$ (see the paragraph following Theorem 2.5). Differential geometric considerations arise naturally in applying Laplace's method to integrate the asymptotic boundary crossing density, and clearly also in the Euler characteristic and tube formulas of excursion sets in [1].

**3. Sums of i.i.d. random variables with multidimensional indices and associated Kolmogorov–Erdős–Feller test.** Let $Y_{\mathbf{k}}, \mathbf{k} \in \mathbf{Z}_+^d$, be i.i.d. random variables with

$$(3.1) \qquad EY_{\mathbf{k}} = 0, \qquad EY_{\mathbf{k}}^2 = 1, \qquad E|Y_{\mathbf{k}}|^3 < \infty.$$

Let $S_{\mathbf{n}} = \sum_{\mathbf{k} \leq \mathbf{n}} Y_{\mathbf{k}}$, where $\mathbf{k} \leq \mathbf{n}$ denotes that $k_i \leq n_i$ for $1 \leq i \leq d$, as in Section 2. Let $|\mathbf{n}| = \prod_{i=1}^n n_i$, $\log \mathbf{n} = (\log n_1, \ldots, \log n_d)$ and $\exp(\mathbf{t}) = (\exp(t_1), \ldots, \exp(t_d))$. Define $X(\log \mathbf{n}) = |\mathbf{n}|^{-1/2} S_{\mathbf{n}}$ and extend the domain of $X$ to $[0,\infty)^d$ by defining $X(\mathbf{t}) = X(\log \mathbf{n})$ when $\log n_i \leq t_i < \log(n_i + 1)$ for all $i$. Let $X_c = X$, $\rho_c = \rho$ and let $D_c$ be a Jordan measurable subset of $\{\mathbf{t} : \sum_i t_i \geq c^3\}$. If $\mathbf{t} = \log \mathbf{n}$ and $\mathbf{t} + \mathbf{u} = \log \mathbf{m}$ for some $\mathbf{m}, \mathbf{n} \in \mathbf{Z}_+^d$, then

$$(3.2) \qquad \begin{aligned} 1 - \rho(\mathbf{t}, \mathbf{t} + \mathbf{u}) &= 1 - \mathrm{Cov}(|\mathbf{n}|^{-1/2} S_{\mathbf{n}}, |\mathbf{m}|^{-1/2} S_{\mathbf{m}}) \\ &= 1 - \exp\left(-\sum_i |u_i|/2\right) \sim \sum_i |u_i|/2 \end{aligned}$$

as $\mathbf{u} \to 0$. From (3.2), it follows that (C) holds with $\alpha = 1$, $L(x) \equiv 1$, $r_t(u) = \sum_i |u_i|/2$, and therefore $\Delta_c = (2c^2)^{-1}$ by (2.5). Moreover, by the Berry–Esseen theorem (cf. [15], Theorem 16.4.1), for $\log \mathbf{n} \leq \mathbf{t} \leq \log(\mathbf{n}+\mathbf{1})$,

$$(3.3) \qquad \left| P\{X(\mathbf{t}) > c - y/c\} - \int_{c-y/c}^{\infty} (2\pi)^{-1/2} e^{-z^2/2} \, dz \right| = O(|\mathbf{n}|^{-1/2})$$

uniformly over $c$ and $y$. Since $\log |\mathbf{n}|/c^2 \to \infty$ uniformly over $\mathbf{t} \in [D_c]_\delta$, it follows from (3.3) that (A1) holds. Moreover, as will be shown in Lemma 3.6, (A3) and (A4) are satisfied uniformly over $[D_c]_\delta$. If we assume in addition that for some $\varepsilon > 0$ and $\kappa > 0$,

$$(3.4) \quad \sup_{\mathbf{t}, \mathbf{u} \in [D_c]_\delta} \|\mathbf{u} - \mathbf{t}\| = O(c^\kappa) \quad \text{and} \quad [D_c]_\delta \subset G_\varepsilon := \left\{\mathbf{t} : t_i \Big/ \sum_j t_j \geq \varepsilon \text{ for all } i\right\},$$

then Lemmas 3.5 and 3.7 show that (A2) and (A5) also hold. Therefore $X(\mathbf{t})$, $\mathbf{t} \in D_c$, is asymptotically Gaussian, and we shall apply Lemma 2.3 and Corollary 2.7 at the end of this section to prove the following two theorems.

THEOREM 3.1. (i) *Assume that for some positive $\zeta_c \to 0$ with $\Delta_c = o(\zeta_c)$,*

$$(3.5) \qquad \zeta_c^d |\{\mathbf{t} \in (\zeta_c \mathbf{Z})^d : I_{t,\zeta_c} \cap \partial D_c \neq \varnothing\}| = o(v(D_c))$$

*as $c \to \infty$. Then $P\{\sup_{\mathbf{t} \in D_c} X(\mathbf{t}) > c\} = O(v(D_c) c^{2d} \psi(c))$.*
(ii) *If (3.4) also holds, then*

$$(3.6) \qquad P\left\{\sup_{\mathbf{t} \in D_c} X(\mathbf{t}) > c\right\} \sim 2^{-d} v(D_c) c^{2d} \psi(c).$$

Theorem 3.1(i) enables us to extend the Kolmogorov–Erdős–Feller test (1.3) to the case of multidimensional time. Let $\beta : \mathbf{Z}_+^d \to (0, \infty)$ be nondecreasing in the sense that $\beta(\mathbf{m}) \leq \beta(\mathbf{n})$ for all $\mathbf{m} \leq \mathbf{n}$. We say that $\beta$ is an upper (lower) class function if

(3.7) $$\sup\{|\mathbf{n}| : |\mathbf{n}|^{-1/2} S_{\mathbf{n}} > \beta(\mathbf{n})\} < (=)\infty \quad \text{a.s.}$$

For $\varepsilon \geq 0$, let $F_\varepsilon = \{\mathbf{n} \in \mathbf{Z}_+^d : \log n_i / \log |\mathbf{n}| \geq \varepsilon \text{ for all } i\}$; in particular $F_0 = \mathbf{Z}_+^d$. Define

(3.8) $$J_\varepsilon = \sum_{\mathbf{n} \in F_\varepsilon} |\mathbf{n}|^{-1} \beta^{2d-1}(\mathbf{n}) e^{-\beta^2(\mathbf{n})/2}.$$

THEOREM 3.2. *If $J_0 < \infty$, then $\beta$ is an upper class function. Conversely, if $J_\varepsilon = \infty$ for some $\varepsilon > 0$, then $\beta$ is a lower class function.*

EXAMPLE 3.3. In the case $d = 1$, since $xe^{-x^2/2}$ is decreasing in $x \geq 1$, it follows that $\int_1^\infty t^{-3/2} b(t) e^{-b^2(t)/2t} \, dt < \infty$ iff $\sum_1^\infty n^{-1} \beta(n) e^{-\beta^2(n)} < \infty$, where $\beta(n) = b(n)/\sqrt{n}$ is nondecreasing. Therefore the integral test (1.3) is equivalent to Theorem 3.2, noting that $J_\varepsilon = J_0$ for all $0 \leq \varepsilon \leq 1$ in the case $d = 1$. Next consider $d = 2$ and let $\beta$ be a positive function on $\mathbf{Z}_+^2$ such that

$$\sum_{n_1 = 1, n_2 \geq 1} n_2^{-1} \beta^3(\mathbf{n}) e^{-\beta^2(\mathbf{n})/2} = \infty,$$

$$\sum_{\mathbf{n}} [n_2^{-1} \beta(\mathbf{n}) \mathbf{1}_{\{n_1 = 1\}} + |\mathbf{n}|^{-1} \beta^3(\mathbf{n}) \mathbf{1}_{\{n_1 \geq 2\}}] e^{-\beta^2(\mathbf{n})/2} < \infty.$$

Then $\sup\{|\mathbf{n}| : |\mathbf{n}|^{-1/2} S_{\mathbf{n}} > \beta(\mathbf{n}), n_1 \geq 2\} < \infty$ a.s. and $\sup\{n_2 : n_2^{-1/2} S_{\mathbf{n}} > \beta(n), n_1 = 1\} < \infty$ a.s., by the first part of Theorem 3.2. On the other hand, $J_0 = \infty$ although $J_\varepsilon < \infty$ for every $\varepsilon > 0$. This shows the importance of using $J_\varepsilon$ instead of $J_0$ for the lower class result in Theorem 3.2.

Let $\beta(\mathbf{n}) = \{(2+\delta) d \log \log |\mathbf{n}|\}^{1/2}$ for $|\mathbf{n}| \geq \mathbf{e}$ and $\delta \geq 0$. Then by the inequality $d^{-1} \sum_{i=1}^d \log n_i \geq (\prod_{i=1}^d \log n_i)^{1/d}$ between arithmetic and geometric means, there exist $C, C' > 0$ such that

$$J_0 \leq C \sum_{\mathbf{n} \in \mathbf{Z}_+^d} (\log \log |\mathbf{n}|)^{d-1/2} / \{|\mathbf{n}| (\log |\mathbf{n}|)^{d(1+\delta/2)}\}$$

$$\leq C' \prod_{i=1}^d \sum_{n_i \in \mathbf{Z}_+} n_i^{-1} (\log n_i)^{-(1+\delta/3)},$$

so $J_0 < \infty$ if $\delta > 0$. Take $0 < \varepsilon < d^{-1}$ and note that the number of $\mathbf{k}$'s such that $\sum_i k_i = m$ and $e^{\mathbf{k}} \in F_\varepsilon$ (so that $k_i \geq \varepsilon m$) is $(B + o(1)) m^{d-1}$ for some $B > 0$. Since $\prod_{i=1}^{d} \sum_{e^{k_i-1} \leq n_i < e^{k_i}} n_i^{-1} \sim 1$ as $\min k_i \to \infty$, it follows that if $\delta = 0$, then there exists $B' > 0$ such that

$$J_\varepsilon \geq \sum_{\mathbf{n} \in F_\varepsilon, \mathbf{n} \geq \mathbf{n}_0} |\mathbf{n}|^{-1}(\log|\mathbf{n}|)^{-d} \geq B' \sum_{m \geq m_0} m^{-1} = \infty.$$

Hence by Theorem 3.2, $\beta(\mathbf{n})$ belongs to the upper class if $\delta > 0$ and to the lower class if $\delta = 0$, yielding the following.

COROLLARY 3.4.  $\limsup_{|\mathbf{n}| \to \infty} S_\mathbf{n}/(2d|\mathbf{n}|\log\log|\mathbf{n}|)^{1/2} = 1$ a.s.

In the case $d = 2$, Zimmerman [38] proved an analogue of Corollary 3.4 for the Brownian sheet, which is a zero-mean Gaussian random field with independent increments and variance function $|\mathbf{t}|$, like that of $S_\mathbf{n}$. His result was subsequently strengthened by Orey and Pruitt ([26], Theorem 2.2) who proved that for the $d$-dimensional Brownian sheet $W(\mathbf{t})$, $P\{W(\mathbf{t})/|\mathbf{t}| \leq f(|\mathbf{t}|)$ for all large $|\mathbf{t}|\} = 1$ (or 0) if

$$(3.9) \qquad \int_1^\infty \xi^{-1}(\log\xi)^{d-1}(\log\log\xi)^{d-1/2} e^{-f^2(\xi)/2} d\xi < \ (\text{or} =)\infty.$$

Actually their result considers $\mathbf{t} \to \mathbf{0}$ rather than $|\mathbf{t}| \to \infty$. However, as $|\mathbf{t}|W(1/t_1, \ldots, 1/t_d)$ is also a Brownian sheet, one can extend their integral test to the preceding statement. Because continuous Gaussian processes (instead of discrete-time sample sums) are involved, the tail distribution of the maximum over a domain $D_c$ does not require condition (3.4); see (2.6) in this connection. Hence unlike (3.8), the integral test (3.9) does not involve $F_\varepsilon$. Instead of the series (3.8), we can rewrite it as an integral when $F_\varepsilon$ is not involved, expressing the convergence criterion in Theorem 3.2 (taking $\varepsilon = 0$) as the integral test

$$(3.10) \qquad \int_1^\infty \cdots \int_1^\infty (t_1 \cdots t_d)^{-1} \beta^{2d-1}(t_1, \ldots, t_d) \\ \times e^{-\beta^2(t_1,\ldots,t_d)/2} dt_1 \cdots dt_d < \ (\text{or} =)\infty.$$

Note that (3.10) considers more general functions $\beta(\mathbf{t})$ than those of the form $f(|\mathbf{t}|)$ considered by Orey and Pruitt [26]. In the case $\beta(\mathbf{t}) = f(|\mathbf{t}|)$, assuming without loss of generality that $c_0 \leq f(\xi)/(\log\log\xi)^{1/2} \leq c_1$ for some $0 < c_0 < c_1$ (see the proof of Theorem 3.2), the change of variables $\xi = t_1 \cdots t_d$ in (3.10) shows that (3.9) and (3.10) are indeed equivalent.

Strong approximations of $S_\mathbf{n}$ have been developed by Rio [32] who has shown that if $Y_\mathbf{k}$, $\mathbf{k} \in \mathbf{Z}_+^d$, are i.i.d. with $EY_\mathbf{k} = 0$, $EY_\mathbf{k}^2 = 1$ and $E|Y_\mathbf{k}|^r < \infty$ for some $r > 2$, then redefining the random variables on a new probability space yields

$$(3.11) \qquad \sup_{0 \leq \mathbf{n} \leq \nu\mathbf{1}} |S_\mathbf{n} - W(\mathbf{n})| = O(\nu^{(d-1)/2}(\log\nu)^{1/2} + \nu^{d/r}) \qquad \text{a.s.}$$

Note that (3.11) bounds the approximation error $S_\mathbf{n} - W(\mathbf{n})$ by

$$\left(\max_{1 \leq i \leq d} n_i\right)^{(d-1)/2} \left(\log \max_{1 \leq i \leq d} n_i\right)^{1/2} + \left(\max_{1 \leq i \leq d} n_i\right)^{d/r},$$

instead of by some sufficiently small power of $|\mathbf{n}| = \prod_{i=1}^d n_i$. Therefore Rio's strong approximation (3.11) cannot be combined with the Orey–Pruitt integral test (3.9) for $W(\mathbf{t})$ to yield a corresponding integral test for $S_\mathbf{n}$. Example 3.3 shows that the integral test (3.9) for $W(\mathbf{t})$ actually does not hold for $S_\mathbf{n}$ which requires a more subtle criterion for a lower class of functions.

The proof of Theorem 3.1 uses the following three lemmas which show that (A2)–(A5) hold under (3.4).

LEMMA 3.5. *Assume* (3.4). *Let* $\mathbf{u}, \mathbf{v} \geq \mathbf{0}$. *Then as* $c \to \infty$,

(3.12)
$$E\{c[X(\mathbf{t} + \mathbf{u}\Delta_c) - X(\mathbf{t})] | X(\mathbf{t}) = c - y/c\}$$
$$\to -\sum_{i=1}^d u_i/4,$$

(3.13)
$$\mathrm{Cov}\{c[X(\mathbf{t} + \mathbf{u}\Delta_c) - X(\mathbf{t})], c[X(\mathbf{t} + \mathbf{v}\Delta_c) - X(\mathbf{t})] | X(\mathbf{t}) = c - y/c\}$$
$$\to \sum_{i=1}^d \min(u_i, v_i)/2,$$

*uniformly over bounded values of* $y$ *and* $\mathbf{t} \in [D_c]_\delta$. *Hence* (A2) *holds*.

PROOF. For $\exp(\mathbf{t}) \in \mathbf{Z}_+^d$ and $\exp(\mathbf{t} + \mathbf{u}\Delta_c) \in \mathbf{Z}_+^d$, define

(3.14)
$$Z_\mathbf{t}(\mathbf{u}) = \sum_{\{\mathbf{k}:\, \mathbf{k} \leq \exp(\mathbf{t}+\mathbf{u}\Delta_c)\} \setminus \{\mathbf{k}:\, \mathbf{k} \leq \exp(\mathbf{t})\}} Y_\mathbf{k}$$
$$= X(\mathbf{t} + \mathbf{u}\Delta_c) \exp\left\{\sum_i (t_i + u_i\Delta_c)/2\right\} - X(\mathbf{t}) \exp\left(\sum_i t_i/2\right).$$

Conditioned on $X(\mathbf{t}) = c - y/c$,

(3.15)
$$c\{X(\mathbf{t} + \mathbf{u}\Delta_c) - X(\mathbf{t})\}$$
$$= cZ_\mathbf{t}(\mathbf{u}) \exp\left\{-\sum_i (t_i + u_i\Delta_c)/2\right\}$$
$$\quad - c(c - y/c)\left\{1 - \exp\left(-\Delta_c \sum_i u_i/2\right)\right\}$$
$$= cZ_\mathbf{t}(\mathbf{u}) \exp\left(-\sum_i t_i/2 - \Delta_c \sum_i u_i/2\right) - \sum_i u_i/4 + o(1).$$

Since $Z_{\mathbf{t}}(\mathbf{u})$ is independent of $X(\mathbf{t})$, (3.12) follows from (3.4) and (3.15). To see this, suppose $\log n_i \leq t_i + u\Delta_c < \log(n_i+1)$ and $\log m_i \leq t_i < \log(m_i+1)$ for $1 \leq i \leq d$. By (3.4), $t_i \geq \varepsilon \sum_j t_j \geq \varepsilon c^3/2$ for $\mathbf{t} \in [D_c]_\delta$ and all large $c$. This implies that $\log(m_i+1) - \log m_i = \log(1+m_i^{-1}) \leq \log(1+e^{-\varepsilon c^3/2}) = o(\Delta_c)$ and that $\log(n_i+1) - \log n_i = o(\Delta_c)$, so by (3.14) and (3.15),

$$E\{c[X(\mathbf{t}+\mathbf{u}\Delta_c) - X(\mathbf{t})]|X(\mathbf{t}) = c - y/c\}$$
$$= -\sum_{i=1}^d (\log n_i - \log m_i)/4\Delta_c + o(1) \to -\sum_{i=1}^d u_i/4.$$

Similarly, for $\mathbf{u}, \mathbf{v} \geq \mathbf{0}$,

$$\text{Cov}(Z_{\mathbf{t}}(\mathbf{u}), Z_{\mathbf{t}}(\mathbf{v})) \sim \prod_{i=1}^d \{\exp(t_i + \min(u_i, v_i)\Delta_c) - \exp(t_i)\}$$
(3.16)
$$\sim \left\{\sum_i \min(u_i, v_i)\Delta_c\right\} \exp\left(\sum_i t_i\right)$$

and (3.13) follows from (3.4), (3.15), (3.16) since $Z_{\mathbf{t}}(\mathbf{v})$ is also independent of $X(\mathbf{t})$. □

LEMMA 3.6. (i) $P\{\max_{\mathbf{k} \leq \mathbf{n}} S_{\mathbf{k}} \geq \lambda\} \leq 2^d P\{S_{\mathbf{n}} \geq \lambda - d(2|\mathbf{n}|)^{1/2}\}$.

(ii) *There exists a positive function $h$, with $\lim_{y\to\infty} h(y) = 0$, satisfying* (A3).

(iii) *There exist nonincreasing functions $N_a$ on $\mathbf{R}_+$ and positive constants $\gamma_a$ such that $\gamma_a \to 0$, $N_a(\gamma_a) + \int_1^\infty y^s N_a(\gamma_a + y)\,dy = o(a^d)$ for all $s > 0$, as $a \to 0$, and* (A4) *holds.*

PROOF. For (i), see [17], Lemma 2.3. To prove (ii), let $\mathbf{u} \geq \mathbf{0}$, $\omega > 0$, $\gamma > 0$. By (3.15),

(3.17) $\quad c[X(\mathbf{t}+\mathbf{u}\Delta_c) - X(\mathbf{t})] \leq (1+o(1))cZ_{\mathbf{t}}(\mathbf{u}) \exp\left(-\sum_i t_i/2\right).$

Since $|\{\mathbf{k}:\mathbf{k} \leq \exp(\mathbf{t}+\mathbf{u}\Delta_c)\}| - |\{\mathbf{k}:\mathbf{k} \leq \exp(\mathbf{t})\}| \sim (\sum_i u_i/2c^2)\exp(\sum_i t_i)$, it follows by the independence of $Z_{\mathbf{t}}(\mathbf{u})$ and $X(\mathbf{t})$ and the Berry–Esseen theorem that for large $c$,

$$P\{c[X(\mathbf{t}+\mathbf{u}\Delta_c) - X(\mathbf{t})] > y' - \gamma|X(\mathbf{t}) = c - y'/c\}$$
(3.18)
$$\leq P\left\{Z_{\mathbf{t}}(\mathbf{u}) > [(y'-\gamma)/2c]\exp\left(\sum_i t_i/2\right)\right\}$$
$$\leq \psi(B(y'-\gamma)) + O\left(c\exp\left(-\sum_i t_i/2\right)\right)$$

for some $B > 0$, uniformly over $\gamma \le y' \le \omega c$. In view of (A1), we can choose $\xi_c \to 0$ such that $P\{X(\mathbf{t}) > c - y'/c\} = (1 + O(\xi_c^2))\psi(c - y'/c)$ uniformly over $\gamma \le y' \le \omega c$. Let $y_j = y + j\xi_c$, $j = 0, 1, \ldots$. Then by (3.18),

$$P\{X(\mathbf{t} + \mathbf{u}\Delta_c) \ge c - \gamma/c, c - \omega \le X(\mathbf{t}) < c - y/c\}$$

$$\le \sum_{0 \le j \le (\omega c - y)/\xi_c} \int_{y_j}^{y_{j+1}} P\{X(\mathbf{t} + \mathbf{u}\Delta_c) > c - \gamma/c | X(\mathbf{t}) = c - y'/c\}$$

$$\times P\{X(\mathbf{t}) \in c - dy'/c\}$$

$$\le \sum_{0 \le j \le (\omega c - y)/\xi_c} \left[ \psi(B(y_j - \gamma)) + O\left( c\exp\left\{ -\sum_j t_j/2 \right\} \right) \right]$$

(3.19)
$$\times [P\{X(\mathbf{t}) > c - y_{j+1}/c\} - P\{X(\mathbf{t}) > c - y_j/c\}]$$

$$\le (1 + o(1)) \sum_{0 \le j \le (\omega c - y)/\xi_c} \left[ \psi(B(y_j - \gamma)) \right.$$

$$\left. + O\left( c\exp\left\{ -\sum_i t_i/2 \right\} \right) \right] \xi_c e^{y_j} \psi(c)$$

$$\le \psi(c) \left[ (1 + o(1)) \right.$$

$$\left. \times \int_y^{\omega c} e^{y'} \psi(B(y' - \gamma))\, dy' + O\left( c\exp\left\{ \omega c - \sum_i t_i/2 \right\} \right) \right].$$

Since $\sum_i t_i \ge c^3$, $c\exp(\omega c - \sum_i t_i/2) = o(1)$. Moreover, $\int_y^\infty e^{y'} \psi(B(y' - \gamma))\, dy' \to 0$ as $y \to \infty$. From (3.17) and (3.18), it follows that for large $c$,

$$P\{[X(\mathbf{t} + \mathbf{u}\Delta_c) - X(\mathbf{t})] > \omega\} \le P\left\{ Z_\mathbf{t}(\mathbf{u}) > (\omega/2) \exp\left( \sum_i t_i/2 \right) \right\}$$

(3.20)
$$\le \psi(Bc\omega) + O\left( c\exp\left( -\sum_i t_i/2 \right) \right)$$

$$= o(\psi(c))$$

if we choose $\omega > B^{-1}$. Hence (ii) follows from (3.19) and (3.20).

To prove (iii), note that $\{\mathbf{k} : \mathbf{k} \le \exp(\mathbf{t} + a\mathbf{1}\Delta_c)\} \setminus \{\mathbf{k} : \mathbf{k} \le \exp(\mathbf{t})\} = \bigcup_{J \subset \{1,\ldots,d\}, J \ne \varnothing} A_J$, where $A_J = \{\mathbf{k} : \exp(t_i) < k_i \le \exp(t_i + a\Delta_c)$ for $i \in J$ and $k_i \le \exp(t_i)$ for $i \notin J\}$. By (i),

$$P\left\{ \sup_{\mathbf{0} \le \mathbf{u} \le a\mathbf{1}} Z_\mathbf{t}(\mathbf{u}) > z \right\}$$

$$(3.21) \quad \leq \sum_{J \subset \{1,\ldots,d\}, J \neq \varnothing} P\left\{\sup_{\mathbf{k} \in A_J} \sum_{\mathbf{m} \leq \mathbf{k}, \mathbf{m} \in A_J} Y_{\mathbf{m}} > z2^{-d}\right\}$$

$$\leq 2^d \sum_{J \subset \{1,\ldots,d\}, J \neq \varnothing} P\left\{\sum_{\mathbf{m} \in A_J} Y_{\mathbf{m}} > z2^{-d} - d(2|A_J|)^{1/2}\right\}.$$

Since $|A_J| \sim (a\Delta_c)^{|J|} \exp(\sum_i t_i)$, it follows by (3.17), (3.21) and the Berry–Esseen theorem [using the same steps as in (3.18)] that for large $c$,

$$P\left\{\sup_{\mathbf{0} \leq \mathbf{u} \leq a\mathbf{1}} c[X(\mathbf{t} + \mathbf{u}\Delta_c) - X(\mathbf{t})] > y' | X(\mathbf{t}) = c - y'/c\right\}$$

$$(3.22) \quad \leq P\left\{\sup_{\mathbf{0} \leq \mathbf{u} \leq a\mathbf{1}} Z_{\mathbf{t}}(\mathbf{u}) > (y'/2c) \exp\left(\sum_i t_i/2\right)\right\}$$

$$\leq 4^d \psi(B'y'/a^{1/2} - 2^{1/2}d) + O\left(c^d \exp\left(-\sum_i t_i/2\right)\right)$$

for some $B' > 0$, uniformly over $\gamma \leq y \leq \omega c$, and therefore

$$P\left\{c - \omega \leq X(\mathbf{t}) < c - \gamma/c, \sup_{\mathbf{0} \leq \mathbf{u} \leq a\mathbf{1}} X(\mathbf{t} + \mathbf{u}\Delta_c) > c\right\}$$

$$(3.23) \quad \leq \psi(c)\left[4^d \int_\gamma^{\omega c} e^{y'} \psi(B'y'/a^{1/2} - 2^{1/2}d)\, dy'\right.$$

$$\left. + O\left(c^d \exp\left(\omega c - \sum_i t_i/2\right)\right)\right].$$

Since $\sum_i t_i \geq c^3$, $O(c^d \exp\{\omega c - \sum_i t_i/2\}) = o(1)$. By (3.17), (3.21), (3.22) and (i),

$$P\left\{\sup_{\mathbf{0} \leq \mathbf{u} \leq a\mathbf{1}} [X(\mathbf{t} + \mathbf{u}\Delta_c) - X(\mathbf{t})] > \omega\right\}$$

$$(3.24) \quad \leq P\left\{\sup_{\mathbf{0} \leq \mathbf{u} \leq a\mathbf{1}} Z_{\mathbf{t}}(\mathbf{u}) > (\omega/2) \exp\left(\sum_i t_i/2\right)\right\}$$

$$\leq 4^d \psi(B'c\omega/a^{1/2} - 2^{1/2}d) + O\left(c^d \exp\left(-\sum_i t_i/2\right)\right)$$

for all large $c$. Let $\gamma_a = a^{1/3}$ and take $\omega > a^{1/2}/B'$ so that $\psi(B'c\omega/a^{1/2} - d2^{1/2}) = o(\psi(c))$. Recalling that $\sum_i t_i \geq c^3$, it follows from (3.23) and (3.24) that for all large $c$ and $\gamma_a \leq \gamma \leq c$,

$$P\left\{\sup_{\mathbf{0} \leq \mathbf{u} \leq a\mathbf{1}} X(\mathbf{t} + \mathbf{u}\Delta_c) > c, X(\mathbf{t}) \leq c - \gamma/c\right\}$$

$$\leq \psi(c) \int_{\gamma}^{\infty} 5^d e^{y'} \psi(B'y'/a^{1/2} - 2^{1/2}d) \, dy' = \psi(c) N_a(\gamma),$$

with $N_a(\gamma_a) + \int_1^{\infty} y^s N_a(\gamma_a + y) \, dy = o(a^p)$ for all $s > 0$ and $p > 0$.  □

LEMMA 3.7. *Assume* (3.4). *For* $\gamma > 0$, *there exist positive constants* $B_1, B_2$ *and* $\eta$ *such that*

(3.25)
$$P\{X(\mathbf{t}) > c - \gamma/c, X(\mathbf{t} + \mathbf{u}\Delta_c) > c - \gamma/c\}$$
$$\leq B_1 \exp(-B_2 \|\mathbf{u}\|^{\eta}) \psi(c - \gamma/c)$$

*uniformly over* $\mathbf{t}, \mathbf{t} + \mathbf{u}\Delta_c \in [D_c]_{\delta}$. *Hence* (A5) *holds*.

PROOF. For $\exp(\mathbf{t}) \in \mathbf{Z}_+^d$ and $\exp(\mathbf{t} + \mathbf{u}\Delta_c) \in \mathbf{Z}_+^d$,

(3.26)
$$P\{X(\mathbf{t}) > c - \gamma/c, X(\mathbf{t} + \mathbf{u}\Delta_c) > c - \gamma/c\}$$
$$\leq P\{X(\mathbf{t}) + X(\mathbf{t} + \mathbf{u}\Delta_c) > 2(c - \gamma/c)\},$$

in which $X(\mathbf{t}) + X(\mathbf{t} + \mathbf{u}\Delta_c)$ is a sum of $\prod_{i=1}^d \max(e^{t_i}, e^{t_i + u_i \Delta_c})$ random variables, each of the form $(\exp(-\sum_i t_i/2) \mathbf{1}_{\{\mathbf{k} \leq \exp(\mathbf{t})\}} + \exp\{-\sum_i (t_i + u_i \Delta_c)/2\} \times \mathbf{1}_{\{\mathbf{k} \leq \exp(\mathbf{t} + \mathbf{u}\Delta_c)\}}) Y_{\mathbf{k}}$. Since $\sum_i t_i \geq c^3$, $X(\mathbf{t}) + X(\mathbf{t} + \mathbf{u}\Delta_c)$ is a sum of at least $\exp(c^3)$ i.i.d. random variables. Using this and $\mathrm{Var}(X(\mathbf{t}) + X(\mathbf{t} + \mathbf{u}\Delta_c)) = 2(1 + \rho(\mathbf{t}, \mathbf{t} + \mathbf{u}\Delta_c))$, we then obtain by the Berry–Esseen theorem that

(3.27)
$$P\{X(\mathbf{t}) + X(\mathbf{t} + \mathbf{u}\Delta_c) > 2(c - \gamma/c)\}$$
$$\leq \psi\left(\left(\frac{2(c - \gamma/c)^2}{1 + \rho(\mathbf{t}, \mathbf{t} + \mathbf{u}\Delta_c)}\right)^{1/2}\right) + O(\exp(-c^3/2)).$$

By (3.2), there exists $\zeta > 0$ such that $1 - \rho(\mathbf{t}, \mathbf{t} + \mathbf{u}\Delta_c) \geq \Delta_c \sum_i |u_i|/4$ if $\Delta_c \sum_i |u_i| \leq \zeta$. Moreover, $1 - \rho(\mathbf{t}, \mathbf{t} + \mathbf{v}) = 1 - \exp(-\sum_i |v_i|/2) \geq \xi := 1 - e^{-\zeta/2}$ if $\sum_i |v_i| \geq \zeta$. Since

$$\psi\left(\left(\frac{2(c - \gamma/c)^2}{1 + \rho(\mathbf{t}, \mathbf{t} + \mathbf{u}\Delta_c)}\right)^{1/2}\right)$$
$$\leq \psi(c - \gamma/c) \left(\frac{1 + \rho(\mathbf{t}, \mathbf{t} + \mathbf{u}\Delta_c)}{2}\right)^{1/2}$$
$$\times \exp\left[\frac{-(c - \gamma/c)^2(1 - \rho(\mathbf{t}, \mathbf{t} + \mathbf{u}\Delta_c))}{2(1 + \rho(\mathbf{t}, \mathbf{t} + \mathbf{u}\Delta_c))}\right]$$
$$\leq \psi(c - \gamma/c) \exp\left[\frac{-(c - \gamma/c)^2(1 - \rho(\mathbf{t}, \mathbf{t} + \mathbf{u}\Delta_c))}{4}\right],$$

it follows from (3.26) and (3.27) that for all large $c$,

$$P\{X(\mathbf{t}) > c - \gamma/c, X(\mathbf{t} + \mathbf{u}\Delta_c) > c - \gamma/c\}$$

(3.28)
$$\leq \psi(c - \gamma/c)$$
$$\times \left[\exp\left(-\sum_i |u_i|/33\right)\mathbf{1}_{\{\Delta_c \sum_i |u_i| \leq \zeta\}} + e^{-c^2 \xi/5}\mathbf{1}_{\{\Delta_c \sum_i |u_i| > \zeta\}}\right].$$

Since $\mathbf{t}$ and $\mathbf{t} + \mathbf{u}\Delta_c$ belong to $[D_c]_\delta$ and since $\sup_{\mathbf{t},\mathbf{v}\in[D_c]_\delta} \|\mathbf{v} - \mathbf{t}\| = O(c^\kappa)$ by (3.4), it follows that $\sum_i |u_i| \leq \sqrt{d}\|u\| = O(c^{\kappa+2})$. Hence (3.25) with $\eta < \min\{1, 2/(\kappa+2)\}$ follows from (3.28). □

PROOF OF THEOREM 3.1. We have already shown that conditions (C), (A1), (A3) and (A4) are satisfied and that (A2) and (A5) also hold under (3.4). Let $\Pi_\mathbf{t} = \{\mathbf{t} + \mathbf{k}\Delta_c \in I_{\mathbf{t},\zeta_c} : \mathbf{k} \in \mathbf{Z}^d\}$. Note that $\zeta_c/\Delta_c \to \infty$ and that

(3.29)
$$P\left\{\sup_{\mathbf{u}\in I_{\mathbf{t},\zeta_c}} X(\mathbf{u}) > c\right\}$$
$$\leq \sum_{\mathbf{u}\in\Pi_\mathbf{t}}\left[P\{X(\mathbf{u}) > c - 1/c\}\right.$$
$$\left.+ P\left\{X(\mathbf{t}) \leq c - 1/c, \sup_{\mathbf{0}\leq\mathbf{v}\leq\mathbf{1}} X(\mathbf{u} + \mathbf{v}\Delta_c) > c\right\}\right]$$
$$= O((\zeta_c/\Delta_c)^d \psi(c)),$$

by (3.3) and (A4). By adding up (3.29) over $\{\mathbf{t} \in (\zeta_c\mathbf{Z})^d : I_{\mathbf{t},\zeta_c} \cap D_c \neq \varnothing\}$, it follows from (3.5) that $P\{\sup_{\mathbf{t}\in D_c} X(\mathbf{t}) > c\} = O(v(D_c)\Delta_c^{-d}\psi(c))$. By (3.2) and Lemma 2.3, $H(t) \equiv 4^{-d}$. If (3.4) also holds, then (A1)–(A5) all hold and Corollary 2.7 can be applied to give (3.6). □

LEMMA 3.8. *Let $\beta:\mathbf{Z}_+ \to (0,\infty)$ be nondecreasing and such that $\beta(\mathbf{n}) \leq \{3d\log\log|\mathbf{n}|\}^{1/2}$. Define $J_\varepsilon$ by (3.8) and let $c(\mathbf{t}) = \beta(\lfloor\exp(\mathbf{t})\rfloor)$ for $\mathbf{t} \in \mathbf{R}_+^d$, $c(\mathbf{0}) = \beta(\mathbf{1})$.*

*(i) If $J_0 < \infty$, then $\sum_{\mathbf{k}\geq\mathbf{0}} c^{2d-1}(\mathbf{k})e^{-c^2(\mathbf{k})/2} < \infty$.*

*(ii) Let $w_1 = 2$ and $w_{j+1} = w_j + \log w_j$ for $j \geq 1$. Then $w_j \sim j\log j$ as $j \to \infty$. For $\mathbf{k} \in \mathbf{Z}_+^d$, define the rectangle $I^{(\mathbf{k})} = \prod_{i=1}^d [w_{k_i}, w_{k_i+1})$ and let $w_\mathbf{k} = (w_{k_1}, \ldots, w_{k_d})$. Assume furthermore that $\beta(\mathbf{n}) \geq \{d\log\log|\mathbf{n}|\}^{1/2}$ and $J_{\varepsilon'} = \infty$ for some $\varepsilon' > 0$. Then for every $0 < \varepsilon < \varepsilon'$,*

(3.30)
$$\sum_{\mathbf{k}\geq\mathbf{3}:\, I^{(\mathbf{k})}\subset G_\varepsilon} v(I^{(\mathbf{k})})c^{2d-1}(w_\mathbf{k})e^{-c^2(w_\mathbf{k})/2} = \infty,$$

*where $G_\varepsilon$ is given in (3.4) and $v(\cdot)$ denotes volume of the rectangle.*

PROOF. Note that $x^{2d-1}e^{-x^2/2}$ is decreasing for $x \geq x_0$. Since $\beta$ is nondecreasing and there are only finitely many $\mathbf{n}$'s with $\beta(\mathbf{n}) < x_0$ in parts (i) and (ii) of the lemma, we can assume without loss of generality that $\beta^{2d-1}(\mathbf{n})e^{-\beta^2(\mathbf{n})/2}$ is decreasing in $\mathbf{n}$. Therefore

$$J_0 \geq \sum_{\mathbf{k} \geq \mathbf{0}} \left\{ \sum_{\mathbf{n}\,:\,e^{k_i} \leq n_i < e^{k_i+1}} |\mathbf{n}|^{-1} \right\} c^{2d-1}(\mathbf{k}+\mathbf{1})e^{-c^2(\mathbf{k}+\mathbf{1})/2}$$

$$\geq \sum_{\mathbf{k} \geq \mathbf{1}} (1-2e^{-1})^d c^{2d-1}(\mathbf{k})e^{-c^2(\mathbf{k})/2},$$

noting that $|\mathbf{n}|^{-1} > \prod_{i=1}^d e^{-(k_i+1)}$ if $e^{k_i} \leq n_i < e^{k_i+1}$ for all $i$ and that $\{\mathbf{n}\,:\,e^{k_i} \leq n_i < e^{k_i+1}$ for all $i\}$ has at least $\prod_{i=1}^d (e^{k_i+1} - e^{k_i} - 1)$ elements. A similar argument also shows that for any $\{i_1, \ldots, i_j\} \subset \{1, \ldots, d\}$,

$$\sum_{\mathbf{k}\,:\,k_{i_1} = \cdots = k_{i_j} = 0, k_i \geq 1 \text{ for } i \notin \{i_1,\ldots,i_j\}} c^{2d-1}(\mathbf{k})e^{-c^2(\mathbf{k})/2} < \infty.$$

To prove (ii), first consider the case $d=1$ for which $J_\varepsilon = J_0$. Here (3.30) follows from $J_0 = \infty$,

$$\sum_{n \geq e^2} n^{-1}\beta(n)e^{-\beta^2(n)/2} \leq \sum_{k=2}^\infty c(w_k)e^{-c^2(w_k)/2} \left\{ \sum_{\exp(w_k) \leq n < \exp(w_{k+1})} n^{-1} \right\},$$

(3.31)
and $\sum_{\exp(w_k) \leq n < \exp(w_{k+1})} n^{-1} \leq w_{k+1} - w_k + 1 = v(I^{(k)}) + 1$.

We next consider the case $d > 1$. Since there are finitely many $\mathbf{n}$'s belonging to $\mathbf{Z}_+^d$ such that $\log \mathbf{n} \in G_{\varepsilon'} \setminus (\bigcup_{\mathbf{k}\,:\,\min k_i \geq 3, I^{(\mathbf{k})} \subset G_\varepsilon} I^{(\mathbf{k})})$, there exists $C > 0$ such that

$$J'_\varepsilon = \sum_{\log \mathbf{n} \in G_{\varepsilon'}} |\mathbf{n}|^{-1}\beta^{2d-1}(\mathbf{n})e^{-\beta^2(\mathbf{n})/2}$$

$$= C + \sum_{\mathbf{k} \geq \mathbf{3}\,:\,I^{(\mathbf{k})} \subset G_\varepsilon} \left\{ \sum_{\log \mathbf{n} \in I^{(\mathbf{k})}} |\mathbf{n}|^{-1}\beta^{2d-1}(\mathbf{n})e^{-\beta^2(\mathbf{n})/2} \right\},$$

which we can bound as in (3.31) to obtain (3.30) if $J_{\varepsilon'} = \infty$, noting that

$$\sum_{\log \mathbf{n} \in I^{(\mathbf{k})}} |\mathbf{n}|^{-1} \leq \prod_{i=1}^d \left( \sum_{\exp(w_{k_i}) \leq n_i < \exp(w_{k_i+1})} n_i^{-1} \right) \leq \prod_{i=1}^d (w_{k_i+1} - w_{k_i} + 1),$$

and that $\prod_{i=1}^d (w_{k_i+1} - w_{k_i} + 1) \sim v(I^{(\mathbf{k})})$ as $\min_{1 \leq i \leq d} k_i \to \infty$. $\square$

PROOF OF THEOREM 3.2. Suppose the theorem holds under the additional assumption $\beta(\mathbf{n}) \leq \{3d \log \log |\mathbf{n}|\}^{1/2}$. To show that it also holds without this additional assumption, define $\hat{\beta}(\mathbf{n}) = \min\{\beta(\mathbf{n}), (3d\log\log|\mathbf{n}|)^{1/2}\}$

for an arbitrary function $\beta: \mathbf{Z}_+^d \to (0, \infty)$. If $J_0(\beta) < \infty$, then $J_0(\hat{\beta}) \leq J_0(\beta) + J_0(\{3d \log \log |\mathbf{n}|\}^{1/2}) < \infty$ and hence

$$\sup\{|\mathbf{n}| : S_\mathbf{n}/|\mathbf{n}|^{1/2} > \beta(\mathbf{n})\} \leq \sup\{|\mathbf{n}| : S_\mathbf{n}/|\mathbf{n}|^{1/2} > \hat{\beta}(\mathbf{n})\} < \infty \quad \text{a.s.}$$

If $J_\lambda(\beta) = \infty$, then $J_\lambda(\hat{\beta}) = \infty$, so $\sup\{|\mathbf{n}| : S_\mathbf{n}/|\mathbf{n}|^{1/2} > \hat{\beta}(\mathbf{n})\} = \infty$ a.s. Since $\sup\{|\mathbf{n}| : S_\mathbf{n}/|\mathbf{n}|^{1/2} > (3d \log \log |\mathbf{n}|)^{1/2}\} < \infty$ a.s., it then follows that $\sup\{|\mathbf{n}| : S_\mathbf{n}/|\mathbf{n}|^{1/2} > \beta(\mathbf{n})\} = \infty$ a.s.

Define $c(\mathbf{t})$ as in Lemma 3.8. In view of the preceding paragraph, we shall assume that $c(\mathbf{t}) \leq \{3d \log(\sum_i t_i)\}^{1/2}$ [and hence $\sum_i t_i \geq c^3(\mathbf{t})$ for large $\mathbf{t}$] and there is no loss of generality. We can apply Theorem 3.1 to $D_c = I_{\mathbf{k},1}$ [noting that (3.5) clearly holds for such cubes with unit width] and combine the result with Lemma 3.8(i) to conclude that if $J_0 < \infty$, then

$$\sum_{\mathbf{k} \geq \mathbf{0}} P\left\{\sup_{\mathbf{t} \in I_{\mathbf{k},1}} X(\mathbf{t}) > c(\mathbf{k})\right\} = O\left(\sum_{\mathbf{k} \geq \mathbf{0}} c^{2d-1}(\mathbf{k}) e^{-c^2(\mathbf{k})/2}\right) < \infty,$$

and therefore $\sum_{\mathbf{k} \geq \mathbf{0}} P\{|\mathbf{n}|^{-1/2} S_\mathbf{n} > \beta(\mathbf{n}) \text{ for some } \mathbf{n} \text{ with } \log \mathbf{n} \in I_{\mathbf{k},1}\} < \infty$. Hence by the Borel–Cantelli lemma, $\beta$ is an upper class function if $J_0 < \infty$.

Suppose $J_\lambda = \infty$ for some $\lambda > 0$. Take $0 < \varepsilon < \lambda$. To prove that $\beta$ is a lower class function, we can assume that $\beta(\mathbf{n}) \geq \{d \log \log |\mathbf{n}|\}^{1/2}$, using an argument similar to that at the beginning of the proof to show that the assumption leads to no loss of generality. For notational simplicity, we focus on the case $d = 2$, as extension of the proof to $d > 2$ is straightforward and the case $d = 1$ does not involve multivariate considerations. Define the rectangles $I^{(\mathbf{k})}$ as in Lemma 3.8(ii) and partition the set $\{\mathbf{k} \geq \mathbf{3} : I^{(\mathbf{k})} \subset G_\varepsilon\}$ of bivariate vectors $\mathbf{k} = (k_1, k_2)$ into four disjoint sets $\mathcal{A}_1, \ldots, \mathcal{A}_4$ so that $k_1$ is odd in $\mathcal{A}_1 \cup \mathcal{A}_2$ (and even in $\mathcal{A}_3 \cup \mathcal{A}_4$) while $k_2$ is odd in $\mathcal{A}_1 \cup \mathcal{A}_3$ (and even in $\mathcal{A}_2 \cup \mathcal{A}_4$). Since the number of $\mathbf{k}$'s belonging to $G_\lambda \setminus (\bigcup_{\mathbf{k} \geq \mathbf{3} : I^{(\mathbf{k})} \subset G_\varepsilon} I^{(\mathbf{k})})$ is finite and since $v(I^{(\mathbf{k})}) = (w_{k_1+1} - w_{k_1})(w_{k_2+1} - w_{k_2}) \sim v(I^{(\mathbf{k-1})})$ as $\min(k_1, k_2) \to \infty$, it follows from Lemma 3.8(ii) that $\sum_{j=1}^4 \sum_{\mathbf{k} \in \mathcal{A}_j} v(I^{(\mathbf{k-1})}) c^{2d-1}(w_\mathbf{k}) e^{-c^2(w_\mathbf{k})/2} = \infty$, and therefore there exists $j$ such that

$$(3.32) \qquad \sum_{\mathbf{k} \in \mathcal{A}_j} v(I^{(\mathbf{k-1})}) c^{2d-1}(w_\mathbf{k}) e^{-c^2(w_\mathbf{k})/2} = \infty.$$

Since $\beta(\mathbf{n}) \geq \{d \log \log |\mathbf{n}|\}^{1/2}$, $c(w_\mathbf{k}) = \beta(\lfloor e^{w_\mathbf{k}} \rfloor) \geq (d + o(1))^{1/2} (\log |w_\mathbf{k}|)^{1/2}$; on the other hand, $w_{k_i} - w_{k_i-1} = \log w_{k_i-1}$, showing that (3.4) holds with $D_c = I^{(\mathbf{k-1})}$, $\kappa = 2$ and $c = c(w_\mathbf{k}) + 2/c(w_\mathbf{k})$. Clearly, $D_c = I^{(\mathbf{k-1})}$ also satisfies (3.5), so Theorem 3.1(ii) can be applied to conclude that

$$\sum_{\mathbf{k} \in \mathcal{A}_j} P\left\{\sup_{\mathbf{t} \in I^{(\mathbf{k-1})}} X(\mathbf{t}) > c(w_\mathbf{k}) + 2/c(w_\mathbf{k})\right\}$$

$$\sim e^{-2} \sum_{\mathbf{k} \in \mathcal{A}_j} v(I^{(\mathbf{k-1})}) c^{2d-1}(w_\mathbf{k}) e^{-c^2(w_\mathbf{k})/2} / (4\sqrt{2\pi}) = \infty,$$

in view of (3.32). This implies that

$$\sum_{\mathbf{k}\in\mathcal{A}_j} P\{S_\mathbf{n}/|\mathbf{n}|^{1/2} > \beta(\mathbf{n}) + 2/c(w_\mathbf{k}) \text{ for some } \mathbf{n} \text{ with } \log\mathbf{n} \in I^{(\mathbf{k-1})}\} \tag{3.33}$$
$$= \infty.$$

Write $\mathbf{n}' < \mathbf{m} \leq \mathbf{n}$ to denote $n'_i < m_i \leq n_i$ for all $i$. For $\log\mathbf{n} \in I^{(\mathbf{k-1})}$, define $S_\mathbf{n}^{(\mathbf{k})} = \sum_{\lfloor\exp(w_{\mathbf{k-2}})\rfloor < \mathbf{m} \leq \mathbf{n}} Y_\mathbf{m}$. We shall show that

$$\sum_{\mathbf{k}\in\mathcal{A}_j} P\left\{\sup_{\log\mathbf{n}\in I^{(\mathbf{k-1})}} |S_\mathbf{n} - S_\mathbf{n}^{(\mathbf{k})}|/|\mathbf{n}|^{1/2} > 1/c(w_\mathbf{k})\right\} < \infty. \tag{3.34}$$

From (3.33) and (3.34), it follows that $\sum_{\mathbf{k}\in\mathcal{A}_j} P(F_\mathbf{k}) = \infty$, where

$$F_\mathbf{k} = \{S_\mathbf{n}^{(\mathbf{k})}/|\mathbf{n}|^{1/2} > \beta(\mathbf{n}) + 1/c(w_\mathbf{k}) \text{ for some } \mathbf{n} \text{ with } \log\mathbf{n} \in I^{(\mathbf{k-1})}\}.$$

Since the $F_\mathbf{k}$ are independent events, $P\{F_\mathbf{k} \text{ i.o.}\} = 1$ by the converse of the Borel–Cantelli lemma. Applying the Borel–Cantelli lemma to (3.34) and combining it with $P\{F_k \text{ i.o.}\} = 1$ then show that $P\{S_\mathbf{n}/|\mathbf{n}|^{1/2} > \beta(\mathbf{n}) \text{ i.o.}\} = 1$. To prove (3.34), let $\mathbf{v} = \lfloor\exp(w_{\mathbf{k-2}})\rfloor$ and note that

$$S_\mathbf{n} - S_\mathbf{n}^{(\mathbf{k})} = \sum_{m_1 \leq n_1, m_2 \leq v_2} Y_\mathbf{m} + \sum_{m_1 \leq v_1, m_2 \leq n_2} Y_\mathbf{m} - \sum_{\mathbf{m}\leq\mathbf{v}} Y_\mathbf{m}.$$

We shall show that

$$\sum_{\mathbf{k}\in\mathcal{A}_j} P\left\{\sup_{\log\mathbf{n}\in I^{(\mathbf{k-1})}} \left|\sum_{m_1\leq n_1, m_2\leq v_2} Y_\mathbf{m}\right|/|\mathbf{n}|^{1/2} > 1/(3c(w_\mathbf{k}))\right\} < \infty. \tag{3.35}$$

Observe that $\sum_{m_1\leq n_1, m_2\leq v_2} Y_\mathbf{m}/|\mathbf{n}|^{1/2} = X(\log n_1, w_{k_2-2})(v_2/n_2)^{1/2}$ and $v_2/n_2 \leq \exp(w_{k_2-2} - w_{k_2-1}) = w_{k_2-2}^{-1}$, and that for large $\mathbf{k}\in\mathcal{A}_j$, $w_{k_2-2}^{1/2}/(3c(w_\mathbf{k})) > (3d\log|w_\mathbf{k}|)^{1/2}$. Therefore

$$P\left\{\sup_{\log\mathbf{n}\in I^{(\mathbf{k-1})}} \left|\sum_{m_1\leq n_1, m_2\leq v_2} Y_\mathbf{m}\right|/|\mathbf{n}|^{1/2} > 1/(3c(w_\mathbf{k}))\right\}$$
$$\leq P\left\{\sup_{w_{\mathbf{k-2}}<\mathbf{t}\leq w_\mathbf{k}} |X(\mathbf{t})| > (3d\log|w_\mathbf{k}|)^{1/2}\right\} = O((k_1+k_2)^{-\lambda})$$

for some $\lambda > 2$, by Theorem 3.1(i). Since $I^{(\mathbf{k-1})} \subset G_\varepsilon$, this proves (3.35). □

**4. Other applications.** Section 2 provides a set of general conditions under which the asymptotic boundary crossing density approximation in Theorem 2.8 is shown to be valid. Given a specific application, one needs only to verify that these assumptions are satisfied. In particular, such verification has been carried out for sums of i.i.d. random variables with multidimensional indices in Section 3, and we begin this section by carrying out similar verification of (C) and (B1)–(B5) for multivariate empirical processes. Let $Y_1, Y_2, \ldots$ be i.i.d. $d$-dimensional random vectors with common distribution function $F$, and let $F_n(t) = n^{-1} \sum_{i=1}^{n} \mathbf{1}_{\{Y_i \leq t\}}$, $t \in \mathbf{R}^d$, be the empirical distribution function of $Y_1, \ldots, Y_n$. Let $Z_n(t) = \sqrt{n}\{F_n(t) - F(t)\}$ be the multivariate empirical process. The limiting distribution of $Z_n$ is that of a Gaussian sheet $Z^0$, for which Adler and Brown [2] proved that

$$(4.1) \quad K_{d,F} c^{2(d-1)} e^{-2c^2} \leq P\Big\{\sup_t Z^0(t) > c\Big\} \leq K_d c^{2(d-1)} e^{-2c^2},$$

where $K_d$ is a constant depending only on $d$ and $K_{d,F}$ is a constant depending on both $d$ and the distribution $F$. For the case $d = 2$ with independent components $Y_{i,1}$ and $Y_{i,2}$ of $Y_i$, $Z^0$ is a pinned Brownian sheet, for which Hogan and Siegmund [20] sharpened (4.1) into

$$(4.2) \quad P\Big\{\sup_t Z^0(t) > c\Big\} \sim (4\log 2) c^2 e^{-2c^2} \quad \text{as } c \to \infty.$$

In Section 5, we apply Corollary 2.9 to prove that if the sample size $n_c$ increases to $\infty$ with $c$ such that $c = o(n_c^{1/6})$, then we can replace $Z^0(t)$ in (4.2) by $Z_{n_c}(t)$ and also extend the result to general $d$ and general distribution $F$ such that

(4.3) $F$ is continuously differentiable and $\partial F/\partial t_i > 0$ for $1 \leq i \leq d$.

In view of (4.3), we can apply a change of variables $t \to F(t)$ and assume that $F$ is a distribution function on the bounded Jordan measurable set $[0,1]^d$, agreeing with the assumptions in Corollary 2.9, whose notation (such as $v_q$) we use in the following theorem.

THEOREM 4.1. *Let* $\mathcal{M} = \{t : F(t) = \frac{1}{2}\}$ *and assume that* (4.3) *holds and* $c = o(n_c^{1/6})$. *Then as* $c \to \infty$,

$$(4.4) \quad P\Big\{\sup_t Z_{n_c}(t) > c\Big\} \\ \sim (8c^2)^{d-1} e^{-2c^2} \int_{\mathcal{M}} \|\nabla F(t)\|^{-1} \prod_{i=1}^{d} \frac{\partial F}{\partial t_i}(t) v_{d-1}(dt),$$

$$(4.5) \quad P\Big\{\sup_t |Z_{n_c}(t)| > c\Big\} \\ \sim 2(8c^2)^{d-1} e^{-2c^2} \int_{\mathcal{M}} \|\nabla F(t)\|^{-1} \prod_{i=1}^{d} \frac{\partial F}{\partial t_i}(t) v_{d-1}(dt).$$

COROLLARY 4.2. (i) *For $d=1$ and continuous distribution function $F$, $P\{\sup_t Z_{n_c}(t) > c\} \sim e^{-2c^2}$ as $c \to \infty$.*

(ii) *For $d=2$, if $F(t_1, t_2) = F_1(t_1)F_2(t_2)$ with continuous univariate distribution functions $F_1$ and $F_2$, then $P\{\sup_t Z_{n_c}(t) > c\} \sim (4\log 2)c^2 e^{-2c^2}$ as $c \to \infty$.*

PROOF. (i) We can assume without loss of generality that $Y_i$ is uniform on $(0,1)$. In this case $\mathcal{M} = \{\frac{1}{2}\}$ and the integral in (4.4) is 1.

(ii) Without loss of generality, assume that $F_1(t_1) = t_1$ and $F_2(t_2) = t_2$, $0 < t_1, t_2 < 1$. In this case $\mathcal{M} = \{(u_1, u_2): u_1 u_2 = \frac{1}{2}, 0 < u_1, u_2 < 1\}$ and the integral in (4.4) becomes

$$\int_{\mathcal{M}} \frac{u_1 u_2}{(u_1^2 + u_2^2)^{1/2}} v_1(du) = \int_{1/2}^1 \frac{1/2}{[u_1^2 + (1/(2u_1))^2]^{1/2}} \left(1 + \frac{1}{4u_1^4}\right)^{1/2} du_1$$

$$= \int_{1/2}^1 \frac{1}{2u_1} du_1 = \frac{1}{2} \log 2,$$

completing the proof for the case $d=2$. □

For $d=1$, Smirnov [33] has shown that $P\{\sup_t Z^0(t) > c\} = e^{-2c^2}$ for all $c > 0$ and Corollary 4.2(i) yields a corresponding asymptotic formula for $Z_{n_c}$, which was used by Chung [12] to prove an upper-lower class theorem for the Kolmogorov–Smirnov statistic. Note that (4.5) says that for some constant $\kappa_d$, $P\{\sup_t |Z_n(t)| > \lambda_n\} \sim \kappa_d \lambda_n^{2(d-1)} e^{-2\lambda_n^2}$ if $\lambda_n \to \infty$ and $\lambda_n = o(n^{1/6})$. Since $n\{F_n(t) - F(t)\}$ is a partial sum of the empirical processes $\mathbf{1}_{\{Y_i \le t\}} - F(t)$, we can apply (4.5) and follow Chung's arguments to prove the following.

COROLLARY 4.3. *Let $\lambda_n$ be a nondecreasing sequence such that $\lambda_n \to \infty$ and let $F$ be a distribution function on $\mathbf{R}^d$ satisfying (4.3). Then $P\{\sup_t |Z_n(t)| > \lambda_n \text{ i.o.}\} = 0$ (or 1) if $\sum_{n=1}^{\infty} n^{-1} \lambda_n^{2d} e^{-2\lambda_n^2} < \infty$ (or $= \infty$).*

An important difference between our proof of the upper-lower class result in Corollary 4.3 and that of Adler and Brown [2] is that they first develop their results for the limiting Kiefer processes and then use the strong approximation theorem of Dudley and Philipp [14] whereas our approach works directly for the empirical process (and of course also for the limiting Kiefer process). The strong approximation approach involves embedding the given process in Brownian motion for which the integral test can be readily shown to hold by using, for example, the tangent approximation (1.1) to the boundary-crossing probability. For partial sums of stationary sequences having long-range dependence, the limiting process is no longer Brownian motion and strong approximation along the lines of Komlós, Major and Tusnády [23], Philipp and Stout [27] and Berkes and Philipp [5] is no longer applicable unless one imposes very restrictive assumptions that are described in the next paragraph. However, the theory in Section 2 can still be applied.

In particular, as in [13] and [19], consider partial sums $S_n := \sum_1^n Y_i$ of linear processes $Y_i := \sum_{j=-\infty}^{\infty} \tau_{i-j}\varepsilon_j$ where $\varepsilon_j$ are i.i.d. random variables with mean 0, variance 1, $Ee^{t|\varepsilon_1|} < \infty$ for some $t > 0$ and $\{\tau_k\}_{k=-\infty}^{\infty}$ satisfies $\sum_{-\infty}^{\infty} \tau_k^2 < \infty$. The sequence $Y_i$ is said to have *long-range dependence* if $E(Y_1 Y_{n+1}) \sim n^{\alpha-2} L(n)$ for some $1 < \alpha < 2$ and slowly varying $L$ so that

(4.6) $$\sigma^2(n) := \text{Var}(S_n) \sim 2\{\alpha(\alpha - 1)\}^{-1} n^\alpha L(n).$$

Defining $Z_n(\cdot)$ by linear interpolation with $Z_n(t) = S_k/\sigma(n)$ for $t = \sigma^2(k)/\sigma^2(n)$, Davydov [13] has shown that $Z_n$ converges weakly to a zero-mean Gaussian process (with correlated increments) whose covariance function is the same as that in (2.1) with $d = 1$ and $r_t \equiv 1$. Although strong approximation theorems are not available for such $S_n$, Chan and Lai [9] have been able to derive integral tests of the type (1.3) for upper-lower class boundaries of $S_n$ in the long-range dependent case by showing that assumptions (C) and (A1)–(A5) of Corollary 2.7 are satisfied by $X_c(t)(= X(t)) = S_{\lfloor e^t \rfloor}/\sigma(\lfloor e^t \rfloor)$ and $D_c = [t_c, t_c^*]$ with $c = o(e^{t_c/6})$, $(t_c^* - t_c)/c^{2/\alpha} \to \infty$ but $t_c^* - t_c = O(c^\kappa)$ for some $\kappa > 0$. Hence application of Corollary 2.7 yields an analog of Theorem 4.1 and therefore also the law of the iterated logarithm (LIL) for partial sums of long-range dependent linear processes; see [9]. Wang, Lin and Gulati [35] recently derived the LIL by using a strong approximation approach that requires $\tau_k$ to have the special form $\tau_k \sim k^{-\beta} L(k)$ as $k \to \infty$ for some $\frac{1}{2} < \beta < 1$ and $\tau_k = 0$ for $k < 0$.

Example 2.2 provides a prototypical example in change-point and signal detection problems, and Theorems 2.4, 2.5, 2.8 and their corollaries can again be applied to a variety of generalizations of Example 2.2 for these applications. Suppose we replace the Brownian motion $W(t)$ by a Gaussian field $X(\mathbf{t})$ and $t_2 - t_1$ is replaced by $\text{Var}(X(\mathbf{t}_2) - X(\mathbf{t}_1))$; see [1] and [7]. Again conditions (C) and (A1)–(A5) can be shown to hold for these applications and also for their discrete-time analogues (like $S_\mathbf{n}$ in Section 3); see [10].

Suppose we replace the Brownian motion $W(t)$ in Example 2.2 by a sample sum process $S_{\lfloor n_c t \rfloor}$, where $S_n = Y_1 + \cdots + Y_n$ and the $Y_i$ are i.i.d. with mean 0, variance 1 and $Ee^{t|Y_1|} < \infty$ for some $t > 0$. Then instead of $X$, we now have a random field $X_c$ defined on $D$ such that

$$X_c(m/n_c, n/n_c) = (S_n - S_m)/(n - m)^{1/2}$$

for $m < n \le an_c$ with $a_1 n_c \le n - m \le a_2 n_c$.

The stopping time $T_c = \inf\{n : \max_{n - a_2 n_c \le m \le n - a_1 n_c} (S_n - S_m)/(n - m)^{1/2} > c\}$ has important applications in sequential change-point detection. Assuming that $n_c/c^6 \to \infty$ as $c \to \infty$ and making use of moderate deviations theory,

Chan and Lai [11] have shown that $X_c$ satisfies conditions (C) and (A1)–(A5). Therefore, analogous to (2.10),

$$(4.7) \qquad P\{T_c \leq an_c\} \sim \psi(c)(c^4/4)[a(a_1^{-1} - a_2^{-1}) - \log(a_2/a_1)].$$

This result provides an important tool for the choice of the threshold $c$ and window sizes of the detection rule $T_c$ to ensure a prescribed false detection rate; see [24] and [11] where the asymptotic optimality (in the sense of quickest detection delay) and extensions (to multivariate $Y_i$ and Markovian $Y_i$) of $T_c$ are also given.

**5. Proofs of Theorems 2.1, 2.4, 2.5, 2.8, 4.1 and their corollaries.** To study the asymptotic distribution (as $t \to \infty$) of $\sup_{0 \leq s \leq t} X(s)$ of a stationary Gaussian process with $EX(s) = 0$, Pickands [28] introduced a method, which has undergone subsequent refinements and is now commonly known as the method of *double sums* (cf. Chapter 2 of [29], [30, 31]), to derive the asymptotic behavior of $P\{\sup_{0 \leq s \leq 1} X(s) > c\}$ as $c \to \infty$. In this section, we modify the double sum method for non-Gaussian fields, to which powerful tools like Slepian's inequality and Fernique's theorem for the Gaussian case (cf. [29]) are no longer applicable. In particular, unlike the traditional double sum $\sum \sum_{i \neq j} P\{\sup_{u \in I^{(i)}} X(u) > c, \sup_{v \in I^{(j)}} X(v) > c\}$ that is shown to be negligible relative to the single sum $\sum_i P\{\sup_{u \in I^{(i)}} X(u) > c\}$ for stationary isotropic Gaussian fields (cf. [28, 29]), note that (2.12) involves $P\{\sup_{u \in I^{(i)}} X_c(u) > c, \sup_{v \in B \setminus I^{(i)}} X_c(v) > c\}$ instead.

The proof of Corollary 2.6 (or 2.7) involves covering $D$ (or $D_c$) by cubes of the form $I_{t,K\Delta_c}$ and using a discrete approximation $A_t(= A_t(a, m, c)) := \{t + \mathbf{k}a\Delta_c : 0 \leq k_i < m, \mathbf{k} \in \mathbf{Z}^d\}$ of $I_{t,K\Delta_c}$, where $a = K/m$. To distinguish from the scalar $K = ma$, we shall use $\mathbf{k}$ to denote the elements of $\mathbf{Z}^d$. Since $P\{\sup_{v \in I_{u,a\Delta_c}} X_c(v) > c, X_c(u) \leq c - \gamma/c\} \leq N_a(\gamma)\psi(c)$ by (A4), approximating the tail probability of $\sup_{u \in I_{t,K\Delta_c}} X_c(u)$ by that of $\sup_{u \in A_t} X_c(u)$ has the error bounds

$$0 \leq \left[P\left\{\sup_{u \in I_{t,K\Delta_c}} X_c(u) > c\right\} - P\left\{\sup_{u \in A_t} X_c(u) > c\right\}\right] \Big/ \psi(c)$$

$$\leq \left[P\left\{c - \gamma/c < \sup_{u \in A_t} X_c(u) \leq c\right\}\right.$$

$$(5.1)$$

$$+ \sum_{u \in A_t} P\left\{\sup_{v \in I_{u,a\Delta_c}} X_c(v) > c, X_c(u) \leq c - \gamma/c\right\}\right] \Big/ \psi(c)$$

$$\leq P\left\{c - \gamma/c < \sup_{u \in A_t} X_c(u) \leq c\right\} \Big/ \psi(c) + (K/a)^d N_a(\gamma),$$

uniformly for $t \in [D]_\delta$ and $\gamma_a \leq \gamma \leq c$. The proof of Theorem 2.4 makes use of (5.1) and Lemma 5.1. Theorem 2.5 is introduced to provide a building block to handle nonstationary random fields (or nonconstant boundaries) in Corollary 2.6 (or Theorem 2.8), which can be proved by much easier arguments in the case of stationary random fields; see Remark 5.1.

LEMMA 5.1. *Under* (C) *and* (A1)–(A3),

$$H_{K,a}(t) := \int_0^\infty e^y P\left\{\sup_{0 \leq k_i < m} W_t(a\mathbf{k}) > y\right\} dy$$

*is uniformly continuous in* $t \in [D]_\delta$ *and* $\sup_{t \in [D]_\delta} H_{K,a}(t) < \infty$. *Moreover, for* $\gamma \geq 0$, *as* $c \to \infty$,

$$\begin{aligned}(5.2)\quad &P\left\{\sup_{u \in A_t} X_c(u) > c - \gamma/c\right\} \\ &\qquad \sim \psi(c - \gamma/c)[1 + H_{K,a}(t)] \qquad \textit{uniformly for } t \in [D]_\delta.\end{aligned}$$

PROOF. Let $\varepsilon > 0$. By (A3), there exists $y^* > \gamma$ such that $h(y^*) < \varepsilon/m^d$ and

$$\begin{aligned}(5.3)\quad 0 &\leq P\left\{\sup_{u \in A_t} X_c(u) > c - \gamma/c\right\} \\ &\quad - P\left\{\sup_{u \in A_t} X_c(u) > c, c - y^*/c < X_c(t) \leq c - \gamma/c\right\} \\ &= P\left\{\sup_{u \in A_t} X_c(u) > c, X_c(t) \leq c - y^*/c\right\} \\ &\leq m^d h(y^*)\psi(c) < \varepsilon \psi(c),\end{aligned}$$

since $|A_t| = m^d$. By (A1), there exists $\xi_c \to 0$ such that

$$(5.4)\qquad |P\{X_c(t) > c - y/c\}/\psi(c - y/c) - 1| = O(\xi_c^2)$$

uniformly for $\gamma \leq y \leq y^*$; we can also assume that $\xi_c^{-1}(y^* - \gamma) \in \mathbf{Z}$. Since $e^{\xi_c} = 1 + \xi_c + O(\xi_c^2)$ and $\psi(c - y/c) \sim e^y \psi(c)$, (5.4) implies that

$$\begin{aligned}(5.5)\quad &P\{c - (y + \xi_c)/c < X_c(t) \leq c - y/c\} \\ &\quad = (1 + O(\xi_c^2))e^{y + \xi_c}\psi(c) - (1 + O(\xi_c^2))e^y \psi(c) \sim \xi_c e^y \psi(c).\end{aligned}$$

By (A2), uniformly for $t \in [D]_\delta$ and $\gamma \leq y \leq y^*$,

$$\begin{aligned}(5.6)\quad &P\left\{\sup_{u \in A_t} X_c(u) > c, c - (y + \xi_c)/c < X_c(t) \leq c - y/c\right\} \\ &\quad \sim P\left\{\sup_{0 \leq k_i < m} W_t(a\mathbf{k}) > y\right\} P\{c - (y + \xi_c)/c < X_c(t) \leq c - y/c\}.\end{aligned}$$

Applying (5.5) to (5.6) and summing (5.6) over $y = j\xi_c + \gamma$ for $j = 0, 1, \ldots,$ $\xi_c^{-1}(y^* - \gamma) - 1$, we obtain (5.2) from (5.3) with arbitrarily small $\varepsilon$. Since $\int_0^\infty e^y P\{W_t(a\mathbf{k}) > y\}\,dy < \infty$ for all $\mathbf{k}$ and $A_t$ is a finite set, $H_{K,a}(t)$ is finite and its uniform continuity follows from (2.1) and (2.2), with the convergence in (2.2) being uniform in $t \in [D]_\delta$ [see the sentence describing assumption (A2)]. Recall in this connection that $\sup_{t \in [D]_\delta, v \in \mathbf{S}^{d-1}} r_t(v) < \infty$, yielding the finiteness of $\sup_{t \in [D]_\delta} H_{K,a}(t)$. □

PROOF OF THEOREM 2.4. Let $a > 0$. By (A4), (5.1) and (5.2), we have for all large $c$,

$$
(5.7) \quad 0 \leq \left[P\left\{\sup_{u \in I_{t,K\Delta_c}} X_c(u) > c\right\} - P\left\{\sup_{u \in A_t} X_c(u) > c\right\}\right]\Big/\psi(c)
$$
$$
\leq 2(e^{\gamma_a} - 1)[1 + H_{K,a}(t)] + (K/a)^d N_a(\gamma_a).
$$

By (A4), for any $\varepsilon > 0$, we can choose $a^*$ small enough such that $N_a(\gamma_a)/a^d < \varepsilon/K^d$ and $2(e^{\gamma_a} - 1) < \varepsilon$ for all $0 < a \leq a^*$. Therefore, by (5.2) and (5.7),

$$
(5.8) \quad (1 - \varepsilon)(1 + H_{K,a}(t))
$$
$$
\leq P\left\{\sup_{u \in I_{t,K\Delta_c}} X_c(u) > c\right\}\Big/\psi(c)
$$
$$
\leq (1 + 2\varepsilon)(1 + H_{K,a}(t)) + \varepsilon,
$$

for all large $c$ and all $t \in [D]_\delta$ and $0 < a \leq a^*$. We shall restrict $a$ and $a^*$ to $\{2^{-j} : j = 1, 2, \ldots\}$ so that the integrand of $H_{K,a}(t)$ is monotone in $a$ and increases to the integrand of $H_K(t)$ as $a \downarrow 0$. Hence by the monotone convergence theorem, $H_{K,a}(t) \to H_K(t)$ as $a \to 0$. Therefore

$$
(5.9) \quad 1 + H_{K,a}(t) \leq 1 + H_K(t) \leq (1 + \varepsilon)(1 + H_{K,a^*}(t)) + \varepsilon,
$$

for all $t \in [D]_\delta$ and $0 < a \leq a^*$ (with $a, a^* \in \{2^{-j} : j = 1, 2, \ldots\}$). We shall use (5.8) and (5.9) in conjunction with Lemma 5.1 to derive the desired conclusions of the theorem.

First note that $M := 1 + \sup_{t \in [D]_\delta} H_K(t) < \infty$ in view of (5.9) and Lemma 5.1 and therefore

$$
(5.10) \quad |H_K(t) - H_{K,a^*}(t)| \leq (M + 1)\varepsilon \quad \text{for all } t \in [D]_\delta,
$$

by (5.9) with $a = a^*$. Because $H_{K,a^*}$ is uniformly continuous by Lemma 5.1,

$$
(5.11) \quad |H_{K,a^*}(t) - H_{K,a^*}(u)| \leq \varepsilon \quad \text{if } \|t - u\| < \delta^*, t, u \in [D]_\delta,
$$

for some $\delta^* > 0$. Since $|H_K(t) - H_K(u)| \leq |H_K(t) - H_{K,a^*}(t)| + |H_K(u) - H_{K,a^*}(u)| + |H_{K,a^*}(t) - H_{K,a^*}(u)|$, it follows from (5.10) and (5.11) that $|H_K(t) - H_K(u)| \leq 2(M+1)\varepsilon + \varepsilon$ if $\|t - u\| < \delta^*$. As $\varepsilon$ is arbitrary, this

shows that $H_K$ is uniformly continuous. Combining (5.8) with (5.10) and the definition of $M$ yields that for all large $c$ and $t \in [D]_\delta$,

$$-\varepsilon - \varepsilon M - \varepsilon^2(M+1) \leq (\psi(c))^{-1} P\Big\{\sup_{u \in I_{t,K\Delta_c}} X_c(u) > c\Big\} - (1 + H_K(t))$$

$$\leq \varepsilon + 2\varepsilon M + 2\varepsilon^2(M+1).$$

Since $\varepsilon$ is arbitrary, this proves Theorem 2.4. $\square$

LEMMA 5.2. *Under* (C) *and* (A1)–(A4), $\sup_{t \in [D]_\delta, K \geq 1} K^{-d} H_K(t) < \infty$ *and* $\{K^{-d} H_K : K \geq 1\}$ *is uniformly equicontinuous on* $[D]_\delta$, *that is,* $\sup_{K \geq 1, t, s \in [D]_\delta, \|t-s\| \leq \varepsilon} |K^{-d} H_K(t) - K^{-d} H_K(s)| \to 0$ *as* $\varepsilon \to 0$.

PROOF. Without loss of generality we can restrict $K$ to be integers. Take any positive integer $a^{-1}$. Note that the integrand of $H_{K,a}(t)$ involves the set $\{a\mathbf{k}: 0 \leq k_i < K/a\}$, which can be partitioned into $K^d$ disjoint subsets $L_j$ such that $|L_j| = a^{-1}$. We can therefore use the arguments at the end of the proof of Lemma 5.1 to bound

$$K^{-d} \sum_{j=1}^{K^d} \Big| P\Big\{\sup_{\mathbf{k} \in L_j} W_t(a\mathbf{k}) > y\Big\} - P\Big\{\sup_{\mathbf{k} \in L_j} W_s(a\mathbf{k}) > y\Big\}\Big|$$

and thereby establish the uniform equicontinuity and boundedness of $\{K^{-d} H_{K,a}: K \geq 1\}$ on $[D]_\delta$. Moreover, by partitioning the cube $[0, K)^d$ similarly into $K^d$ unit cubes, it can be shown that $\sup_{K \geq 1, t \in [D]_\delta} |K^{-d} H_K(t) - K^{-d} H_{K,a}(t)| \to 0$ as $a \to 0$. Hence we can proceed as in (5.10) and (5.11) but with $H_{K,a}$ and $H_K$ replaced by $K^{-d} H_{K,a}$ and $K^{-d} H_K$ to prove the uniform equicontinuity and boundedness of $\{K^{-d} H_K : K \geq 1\}$. $\square$

LEMMA 5.3. *Under* (C) *and* (A1)–(A5), *there exist constants* $s_K \to 0$ *as* $K \to \infty$ *such that*

$$(5.12) \quad P\Big\{\sup_{u \in I_{t,K\Delta_c}} X_c(u) > c, \sup_{v \in B \setminus I_{t,K\Delta_c}} X_c(v) > c\Big\} \leq s_K K^d \psi(c)$$

*for $c$ large enough, uniformly over $t \in [D]_\delta$ and over subsets $B$ of $[D]_\delta$ with bounded volume.*

PROOF. Let $a > 0$ and $0 < q < p$. Then

$$G_a := \sum_{w \in (a\mathbf{Z})^d} \exp(\|w\|^q) f(\|w\|) < \infty.$$

Let $m, n$ be positive integers that are large enough such that

$$\sum_{w \in (a\mathbf{Z})^d, \|w\| \geq na} \exp(\|w\|^q) f(\|w\|) < \varepsilon a^d/8$$

and
$$[1 - (1 - 2n/m)^d] < \varepsilon a^d/8G_a.$$

Let $K = ma$, $F_{1,t} = \{t + \mathbf{k}a\Delta_c : n \leq k_i < m - n, \mathbf{k} \in \mathbf{Z}^d\}$, $F_{2,t} = A_t \setminus F_{1,t}$, $B_t = \{t + a\mathbf{k}\Delta_c \in B \setminus I_{t,K\Delta_c} : \mathbf{k} \in \mathbf{Z}^d\}$, $g_{uv} = \min\{c - \gamma_a, (\|v-u\|/\Delta_c)^q\}$. Then by (A5),

$$P\{X_c(u) > c - (\gamma_a + g_{u,v,c})/c, X_c(v) > c - (\gamma_a + g_{u,v,c})/c\}$$
(5.13)
$$\leq \psi(c - (\gamma_a + g_{uv})/c)f(\|u-v\|/\Delta_c)$$
$$\leq 2e^{g_{uv}}\psi(c - \gamma_a/c)f(\|u-v\|/\Delta_c),$$

for all large $c$ and small $a$. For $u \in F_{1,t}$ and $v \in B_t$, $\|u-v\|/\Delta_c \geq na$ and $g_{uv} \leq (\|u-v\|/\Delta_c)^q$. Noting that $|F_{1,t}| \leq m^d$, $|F_{2,t}| \leq m^d - (m-2n)^d = m^d[1 - (1-2n/m)^d]$, and that $\sum_{u \in A_t} = \sum_{j=1}^{2}\sum_{u \in F_{j,t}}$, we obtain from (5.13) that for all large $c$ and small $a$,

$$\sum_{u \in A_t}\sum_{v \in B_t} P\{X_c(u) > c - (\gamma_a + g_{uv})/c, X_c(v) > c - (\gamma_a + g_{uv})/c\}$$

$$\leq \psi(c - \gamma_a/c)m^d$$
(5.14)
$$\times \left\{\sum_{w \in (a\mathbf{Z})^d, \|w\| \geq na} \exp(\|w\|^q)f(\|w\|) + [1 - (1-2n/m)^d]G_a\right\}$$

$$< (\varepsilon K^d/2)\psi(c - \gamma_a/c).$$

Define $\lambda_w = \min_{u \in A_t} g_{uw}$ if $w \in B_t$, and $\lambda_w = 0$ if $w \in A_t$. Then

$$P\left\{\sup_{u \in I_{t,K\Delta_c}} X_c(u) > c, \sup_{v \in B \setminus I_{t,K\Delta_c}} X_c(v) > c\right\}$$

(5.15) $$\leq \sum_{u \in A_t}\sum_{v \in B_t} P\{X_c(u) > c - (\gamma_a + g_{uv})/c, X_c(v) > c - (\gamma_a + g_{uv})/c\}$$

$$+ \sum_{w \in A_t \cup B_t} P\left\{\sup_{z \in I_{w,a\Delta_c}} X_c(z) > c, X_c(w) \leq c - (\gamma_a + \lambda_w)/c\right\}.$$

On the right-hand side of (5.15), the first sum can be bounded by (5.14) and the second sum by

$$\sum_{u \in A_t} P\left\{\sup_{z \in I_{u,a\Delta_c}} X_c(z) > c, X_c(u) \leq c - \gamma_a/c\right\}$$

(5.16) $$+ \sum_{v \in B_t} P\left\{\sup_{z \in I_{v,a\Delta_c}} X_c(z) > c, X_c(v) \leq c - (\gamma_a + \lambda_v)/c\right\}$$

$$\leq (K/a)^d N_a(\gamma_a)\psi(c) + \sum_{v \in B_t} N_a(\gamma_a + \lambda_v)\psi(c),$$

in view of (A4) and that $|A_t| = m^d = (K/a)^d$. To bound the last sum $\sum_{u \in B_t}$ in (5.16), first consider the case $d = 1$. Since $\lambda_v \geq \min\{(ak)^q, c - \gamma_a\}$ if $ak\Delta_c \leq \inf_{u \in A_t} |v - u| < a(k+1)\Delta_c$, and since $N_a$ is nonincreasing, it follows that

$$\sum_{v \in B_t} N_a(\gamma_a + \lambda_v)\psi(c)$$

(5.17)
$$\leq 2\left\{\sum_{k=1}^{\infty} N_a(\gamma_a + (ak)^q) + N_a(c)v(B)/(a\Delta_c)\right\}$$

$$\leq 2\left\{a^{-1}\int_0^{\infty} N_a(\gamma_a + y^q)\,dy + v(B)N_a(c)/(a\Delta_c)\right\}.$$

Integration by parts shows that the integral in (5.17) approaches 0 as $a \to 0$, since $N_a(\gamma_a) + \int_1^{\infty} w^s N_a(\gamma_a + w)\,dw = o(a)$ for $s > (q^{-1} - 1)^+$. Moreover, in view of (2.5), $N_a(c)/\Delta_c = O(\int_{c/2-\gamma_a}^{c-\gamma_a} w^s N_a(\gamma_a + w)\,dw) = o(a)$ as $a \to 0$ and $c \to \infty$, for $s > 2/\alpha$. Therefore, $\sum_{v \in B_t} N_a(\gamma_a + \lambda_v) \leq \varepsilon/4$ for all large $c$ and small $a$. In general, for $d > 1$,

$$\sum_{v \in B_t} N_a(\gamma_a + \lambda_v)$$

(5.18)
$$\leq 2\left\{d\sum_{j=1}^{\infty}(a^{-1}K + 2j)^{d-1}N_a(\gamma_a + (aj)^q) + N_a(c)v(B)/(a\Delta_c)^d\right\}$$

$$\leq \varepsilon K^{d-1}/4 \qquad \text{for all large } c \text{ and small } a,$$

as can be shown by arguments similar to those in the case $d = 1$. Combining (5.15) with (5.14) and (5.16)–(5.18) yields the desired conclusion. $\square$

PROOF OF THEOREM 2.5. Let $1 > \varepsilon > 0$. There exists $K^*$ such that $s_K \leq \varepsilon/3$ for all $K \geq K^*$. For fixed $t \in D$ and $K \geq K^*$, define

(5.19)
$$\underline{\Lambda} = \{u \in (K\Delta_c\mathbf{Z})^d : I_{u,K\Delta_c} \subset I_{t,\ell_c\Delta_c}\},$$
$$\overline{\Lambda} = \{u \in (K\Delta_c\mathbf{Z})^d : I_{u,K\Delta_c} \cap I_{t,\ell_c\Delta_c} \neq \varnothing\}, \qquad J_u = I_{u,K\Delta_c}.$$

Covering $I_{t,\ell_c\Delta_c}$ by cubes of length $K\Delta_c$ and letting $B$ be a subset of $[D]_\delta$ containing $I_{t,\ell_c\Delta_c}$ and such that $v(B) \leq v_0$, we have

(5.20)
$$\sum_{u \in \underline{\Lambda}}\left[P\left\{\sup_{v \in J_u} X_c(v) > c\right\} - P\left\{\sup_{v \in J_u} X_c(v) > c, \sup_{w \in B\setminus J_u} X_c(w) > c\right\}\right]$$

$$\leq P\left\{\sup_{v \in I_{t,\ell_c\Delta_c}} X_c(v) > c\right\} \leq \sum_{u \in \overline{\Lambda}} P\left\{\sup_{v \in J_u} X_c(v) > c\right\}.$$

By Theorem 2.4 and Lemma 5.3, as $c \to \infty$,

$$(1 + o(1))\psi(c) \sum_{u \in \underline{\Lambda}} [H_K(u) - s_K K^d]$$

(5.21)
$$\leq P\left\{\sup_{v \in I_{t,\ell_c \Delta_c}} X_c(v) > c\right\} \leq (1 + o(1))\psi(c) \sum_{u \in \overline{\Lambda}} H_K(u),$$

uniformly in $t \in D$. In view of $\ell_c \Delta_c \to 0$ and the uniform equicontinuity in Lemma 5.2, we can choose $c^*$ large enough so that $|K^{-d}H_K(u) - K^{-d}H_K(t)| \leq \varepsilon/3$ for all $c \geq c^*$, $\sqrt{\ell_c} \geq K \geq K^*$, $t \in D$ and $u \in \overline{\Lambda}$ ($= \overline{\Lambda}(t; K\Delta_c, \ell_c \Delta_c)$). Putting this and the bound $s_K \leq \varepsilon/3$ in (5.21) and dividing (5.21) by $\ell_c^d \psi(c)$, we obtain for all $c \geq c^*$, $\sqrt{\ell_c} \geq K \geq K^*$ and $t \in D$,

$$(1 - \varepsilon)\{K^{-d}H_K(t) - 2\varepsilon/3\} \leq P\left\{\sup_{v \in I_{t,\ell_c \Delta_c}} X_c(v) > c\right\} / (\ell_c^d \psi(c))$$

(5.22)
$$\leq (1 + \varepsilon)\{K^{-d}H_K(t) + \varepsilon/3\},$$

since $|\overline{\Lambda}| \sim |\underline{\Lambda}| \sim (\ell_c/K)^d$. By Lemma 5.2, $M := \sup_{t \in [D]_\delta, K \geq 1} K^{-d}H_K(t) < \infty$. Therefore, it follows from (5.22) that

$$(5.23) \sup_{t \in D} \left| P\left\{\sup_{v \in I_{t,\ell_c \Delta_c}} X_c(v) > c\right\} / (\ell_c^d \psi(c)) - K^{-d}H_K(t) \right| \leq \varepsilon M + 2\varepsilon/3,$$

for all $c \geq c^*$ and $\sqrt{\ell_c} \geq K \geq K^*$. Letting $c \to \infty$ in (5.23) yields

$$\sup_{t \in D} |K^{-d}H_K(t) - \tilde{K}^{-d}H_{\tilde{K}}(t)| \leq 2\varepsilon M + 4\varepsilon/3,$$

if $K \geq K^*$ and $\tilde{K} \geq K^*$, establishing that $\{K^{-d}H_K\}$ is uniformly Cauchy. Hence $K^{-d}H_K(t)$ converges uniformly in $t \in D$ to $H(t)$, which is also bounded by $M$. We can therefore proceed as in the second paragraph of the proof of Theorem 2.4 to show that $H(t)$ is uniformly continuous in $t \in D$. Moreover, taking $K$ large enough such that $\sup_{t \in D} |K^{-d}H_K(t) - H(t)| \leq \varepsilon/3$, it follows from (5.23) that

$$\sup_{t \in D} \left| P\left\{\sup_{v \in I_{t,\ell_c \Delta_c}} X_c(v) > c\right\} / (\ell_c^d \psi(c)) - H(t) \right| \leq \varepsilon(M + 1)$$

for all $c \geq c^*$, proving (2.11).

We next show that $\inf_{t \in D} H(t) > 0$. For the function $f$ in (A5), we can choose $a > 0$ large enough so that $\sum_{\mathbf{k} \neq \mathbf{0}} f(a\mathbf{k}) \leq 1/2$. Let $K = ma$ and define $A_t = A_t(a, m, c)$ as in the paragraph preceding Lemma 5.1 so that $|A_t| = m^d$. Then by (A1) and (A5), as $c \to \infty$,

$$P\left\{\sup_{u \in A_t} X_c(u) > c\right\}$$

$$
\begin{aligned}
&\geq \sum_{u \in A_t} \left[ P\{X_c(u) > c\} - \sum_{v \in A_t, v \neq u} P\{X_c(u) > c, X_c(v) > c\} \right] \\
&\geq \sum_{u \in A_t} (1 + o(1)) \psi(c)/2 \\
&= (1 + o(1)) m^d \psi(c)/2,
\end{aligned}
$$
(5.24)

uniformly in $t \in D$ and $m \geq 2$. Combining (5.24) with Theorem 2.4 yields

$$
1 + H_{ma}(t) = \lim_{c \to \infty} (\psi(c))^{-1} P \left\{ \sup_{u \in I_{t,ma\Delta_c}} X_c(u) > c \right\}
$$

$$
\geq \limsup_{c \to \infty} (\psi(c))^{-1} P \left\{ \sup_{u \in A_t} X_c(u) > c \right\} \geq m^d/2
$$

for all $m \geq 2$ and $t \in D$. Since $\lim_{K \to \infty} K^{-d} H_K(t) = H(t)$, it then follows that $H(t) \geq a^{-d}/2$ for all $t \in D$. Finally, to prove (2.12), apply (5.12) to obtain that for all $t \in D$ and large $c$,

$$
P \left\{ \sup_{u \in I_{t,\ell_c \Delta_c}} X_c(u) > c, \sup_{v \in B \setminus I_{t,\ell_c \Delta_c}} X_c(v) > c \right\}
$$

$$
\leq \sum_{u \in \overline{\Lambda}} P \left\{ \sup_{v \in J_u} X_c(v) > c, \sup_{v \in B \setminus J_u} X_c(v) > c \right\}
$$

$$
\leq |\overline{\Lambda}| s_K K^d \psi(c).
$$

Since $s_K \to 0$ as $K \to \infty$ and $|\overline{\Lambda}| \sim (\ell_c/K)^d$ as $\ell_c/K \to \infty$, (2.12) follows. □

PROOF OF COROLLARY 2.6. A basic idea of the proof is to cover the bounded, Jordan measurable set $D$ by cubes of length $\ell_c \Delta_c$, with $\ell_c \to \infty$ such that $\ell_c \Delta_c \to 0$. Define $\underline{\Lambda}$, $\overline{\Lambda}$ and $J_u$ as in (5.19) but with $(K \Delta_c \mathbf{Z})^d$ replaced by $(\ell_c \Delta_c \mathbf{Z})^d$, $I_{u,K\Delta_c}$ by $I_{u,\ell_c \Delta_c}$, and $I_{t,\ell_c \Delta_c}$ by $D$. Then (5.20) still holds with these new definitions of $\underline{\Lambda}$, $\overline{\Lambda}$ and $J_u$ and also with $B$ replaced by the bounded set $[D]_\delta$. Labeling it as (5.20′), the upper and lower bounds in (5.20′) are both asymptotically equivalent to $(\ell_c \Delta_c)^{-d} \ell_c^d \psi(c) \int_D H(t)\, dt$ by Theorem 2.5, since $\ell_c \Delta_c \to 0$ and $H(t)$ is continuous. □

REMARK 5.1. Corollary 2.6 can be proved by easier arguments, to be sketched below, when $X_c(t) = X(t)$ is stationary. Let

$$
\Lambda_J = \{t \in D : t \in (J\Delta_c)^d\},
$$

$$
F = \left\{ \sup_{u \in D} X(u) > c \right\},
$$

$$F_t = \left\{ \sup_{u \in A_t(a, K/a, c)} X(u) > c - \gamma/c \right\},$$

$$\widetilde{F}_t = \left\{ \sup_{u \in A_t(\tilde{a}, \widetilde{K}/\tilde{a}, c)} X(u) > c - \gamma/c \right\}.$$

Then

(5.25)
$$\left| \sum_{t \in \Lambda_K} P(F_t) - \sum_{t \in \Lambda_{\widetilde{K}}} P(\widetilde{F}_t) \right|$$
$$\leq \left| \sum_{t \in \Lambda_K} P(F_t) - P(F) \right| + \left| \sum_{t \in \Lambda_{\widetilde{K}}} P(\widetilde{F}_t) - P(F) \right|.$$

It can be shown by arguments similar to those in the proof of Lemmas 2.3 and 2.4 of [31] that $\limsup_{c \to \infty} |\sum_{t \in \Lambda_K} P(F_t) - P(F)|/\{\psi(c)\Delta_c^{-d}\} \to 0$ as $K \to \infty$, $a \to 0$ and $\gamma \to 0$. Moreover, by Lemma 5.1 and stationarity,

(5.26) $$\sum_{t \in \Lambda_K} P(F_t) \sim v(D)(K\Delta_c)^{-d}\psi(c - \gamma/c)(1 + H_{K,a}),$$

and a similar relation also holds for $\sum_{t \in \Lambda_{\widetilde{K}}} P(\widetilde{F}_t)$. Hence by (5.25),

(5.27) $$|K^{-d} H_{K,a} - \widetilde{K}^{-d} H_{\widetilde{K},\tilde{a}}| \to 0 \quad \text{as } K, \widetilde{K} \to \infty, a, \tilde{a} \to 0,$$

which implies that $\lim_{K \to \infty, a \to 0} K^{-d} H_{K,a}$ exists by the Cauchy convergence property, yielding $H$ as the limit. For nonstationary random fields, we do not have the simple relation (5.26) and cannot show the existence of the limit of $K^{-d} H_{K,a}(t)$ via Cauchy convergence as in (5.27). This is why more complicated arguments are needed in the proofs of Lemma 5.3 and Theorem 2.5, from which Corollary 2.6 follows. Concerning the proofs of Theorems 2.4 and 2.5, since $H_{K,a}$ and $H$ are defined by the Gaussian processes $W_t$ rather than the process $X_c$ satisfying (A1)–(A5), one may wonder why these assumptions have been involved in their proofs (and also that of Lemma 5.1) to establish continuity and boundedness properties of $H_{K,a}$ and $H$. It turns out that for a Gaussian random field $X_c = X$ whose covariance function satisfies condition (C), assumptions (A1)–(A5) also hold with $W_t$ being the limiting process in (A2); see the following proof of Theorem 2.1 which generalizes the Qualls–Watanabe result (2.6) to nonstationary Gaussian fields.

PROOF OF COROLLARY 2.7. Here we modify (5.19) into

(5.28)
$$\underline{\Lambda}_c = \{u \in (\zeta_c \mathbf{Z})^d : I_{u,\zeta_c} \subset D_c\},$$
$$\overline{\Lambda}_c = \{u \in (\zeta_c \mathbf{Z})^d : I_{u,\zeta_c} \cap D_c = \varnothing\},$$
$$J_u = I_{u,\zeta_c},$$

and replace $B$ in (5.20) by $[D_c]_\delta$ so that we have here a corresponding version of (5.20), labeled as (5.20*). Apply (2.12) with $\ell_c = \zeta_c/\Delta_c$ together with (2.14) and (5.20*) to derive (2.15). □

PROOF OF THEOREM 2.8. It follows from (2.11) and (2.16) that

$$(5.29) \quad P\{X_c(u) > b_c(u) \text{ for some } u \in I_{t,\zeta_c}\} \sim (\zeta_c/\Delta_{b_c(t)})^d \psi(b_c(t)) H(t).$$

From (2.12) with the boundary $c$ replaced by $\underline{b}_c(t)$ and with $\ell_c = \zeta_c/\Delta_c$, it follows that

$$(5.30) \quad \begin{aligned} &P\{X_c(u) > b_c(u), X_c(v) > \max(b_c(v), \underline{b}_c(t)) \\ &\qquad \text{for some } u \in I_{t,\zeta_c}, v \in B \setminus I_{t,\zeta_c}\} \\ &\leq P\left\{\sup_{u \in I_{t,\zeta_c}} X_c(u) > \underline{b}_c(t), \sup_{v \in B \setminus I_{t,\zeta_c}} X_c(v) > \underline{b}_c(t)\right\} \\ &= o((\zeta_c/\Delta_{b_c(t)})^d \psi(b_c(t))) \end{aligned}$$

uniformly over $t \in D_c$ and over subsets $B$ of $[D]_\delta$ with bounded volume. Then

$$(5.31) \quad \begin{aligned} \sum_{u \in \underline{\Lambda}_c} &\left[ P\{X_c(w) > b_c(w) \text{ for some } w \in J_u\} \right. \\ &\quad - \sum_{v \in \underline{\Lambda}_c, v \neq u} P\{X_c(w) > b_c(w), X_c(z) > \max(b_c(z), \underline{b}_c(u)) \\ &\qquad\qquad\qquad\qquad\qquad\qquad \text{for some } w \in J_u, z \in J_v\} \Bigg] \\ &\leq P\{X_c(u) > b_c(u) \text{ for some } u \in D_c\} \\ &\leq \sum_{u \in \overline{\Lambda}_c} P\{X_c(w) > b_c(w) \text{ for some } w \in J_u\}, \end{aligned}$$

in which $\underline{\Lambda}_c$ and $\overline{\Lambda}_c$ are defined by (5.28). By (5.29) and (5.30), the lower and upper bounds in (5.31) are asymptotically equal. Since $\zeta_c \to 0$, the desired conclusion then follows. □

PROOF OF COROLLARY 2.9. Since $\alpha < 2$, there exists $\varepsilon > 0$ such that $\Delta_c = o(c^{-(1+\varepsilon)})$. We can therefore choose $\zeta_c \to 0$ and $\xi_c \to 0$ such that

$$(5.32) \quad \zeta_c/\Delta_c \to \infty, \qquad \xi_c \geq c^{-1} \log c, \qquad \zeta_c \xi_c = o(c^{-2}),$$

so $\zeta_c = o(\xi_c)$. Consider the tubular neighborhood $U_{\xi_c}$ of $\mathcal{M}$. For sufficiently small $\xi$, the elements of $U_\xi$ can be uniquely represented in the form $x + y$ with $x \in \mathcal{M}$, $y \in T\mathcal{M}^\perp(y)$ and $\|y\| < \xi$. Since $\nabla b(t) = 0$ for all $t$ belonging to the compact set $\mathcal{M}$, there exists $B > 0$ such that $\|\nabla b(u)\| \leq B\xi_c$ for all $u \in [U_{\xi_c}]_{2\zeta_c}$. Combining this with (5.32) yields

$$\sup_{u \in I_{t,\zeta_c}} b(u) - \inf_{u \in I_{t,\zeta_c}} b(u) = O(\zeta_c \xi_c) = o(c^{-2})$$

(5.33)

uniformly over $t \in [U_{\xi_c}]_{2\zeta_c}$.

Recalling that $b_c(u) = cb(u)$ and applying the identity $y^2 - x^2 = (y-x)(y+x)$, we can conclude from (5.33) that $\sup_{t \in [U_{\xi_c}]_{2\zeta_c}} [\overline{b}_c^2(t) - \underline{b}_c^2(t)] = o(1)$.

Let $y \in \mathcal{M}$ and $z \in T\mathcal{M}^\perp(y)$. Then $b(y+z) = b_D + z'\nabla^2 b(y)z/2 + O(\|z\|^3)$. Applying Theorem 2.8 to $D_c = U_{\xi_c}$ yields

$$P\{X_c(t) > cb(t) \text{ for some } t \in U_{\xi_c}\}$$

(5.34)

$$\sim \Delta_c^{-d} \int_{U_{\xi_c}} \psi(cb(t))(b(t))^{2d/\alpha} H(t)\, dt$$

$$\sim \Delta_c^{-d} \psi(cb_D) b_D^{2d/\alpha}$$

$$\times \int_\mathcal{M} H(y) \int_{z \in T\mathcal{M}^\perp(y), \|z\| \leq \xi_c} \exp(-c^2 b_D z' \nabla^2 b(y)z/2)\, dz v_q(dy).$$

Since $\nabla_\perp^2 b(y)$ is positive definite, $\inf_{u \in D \setminus U_{\xi_c}} b(u) \geq b_D + B'\xi_c^2$ for some $B' > 0$. Hence by Theorem 2.8,

$$P\{X_c(t) > cb(t) \text{ for some } t \in D \setminus U_{\xi_c}\}$$

(5.35)

$$\leq P\left\{\sup_{t \in D \setminus U_{\xi_c}} X_c(t) > c(b_D + B'\xi_c^2)\right\}$$

$$= o(\Delta_c^{-d} \psi(cb_D)),$$

in view of (5.32). Combining (5.34) with (5.35) and evaluating the inner integral in (5.34) give (2.19). □

PROOF OF THEOREM 2.1. We shall show that (A1)–(A5) hold for the Gaussian random field $X$ satisfying (2.7), which is the same as condition (C) in the present case $X_c = X$, and hence Theorem 2.1 follows from Theorem 2.5 and Corollary 2.6. In particular, (A1) follows from the well-known asymptotic tail behavior of a normal distribution. Let $\rho(t,u) = E[X(t)X(u)]$. Since the conditional distribution of $X(t + u\Delta_c)$ given $X(t)$ is normal with mean $\rho(t, t + u\Delta_c)X(t)$, it follows from (2.7) that as $c \to \infty$,

$$E\{c[X(t + u\Delta_c) - X(t)] | X(t) = c - y/c\}$$

(5.36)

$$= -c[1 - \rho(t, t + u\Delta_c)](c - y/c)$$

$$\to -\|u\|^\alpha r_t(u/\|u\|)/2,$$

$$\operatorname{Cov}\{c[X(t+u\Delta_c) - X(t)], c[X(t+v\Delta_c) - X(t)] | X(t) = c - y/c\}$$

(5.37)
$$= c^2[\rho(t+u\Delta_c, t+v\Delta_c) - \rho(t, t+v\Delta_c)\rho(t, t+u\Delta_c)]$$

$$\to [-\|v-u\|^\alpha r_t((v-u)/\|v-u\|)$$

$$+ \|v\|^\alpha r_t(v/\|v\|) + \|u\|^\alpha r_t(u/\|u\|)]/2.$$

Since $\{c[X(t+ak\Delta_c) - X(t)] : 0 \leq k_i < m\}$ is multivariate normal, (A2) then follows. Let $\gamma > 0$. Since $\psi(c - z/c) \sim e^z \psi(c)$ for all $z \geq 0$ and there exist constants $B, B' > 0$ such that $P\{W_t(u) > z - \gamma\} \leq B \exp(-B' z^2)$, it follows from (5.36) and (5.37) that as $c \to \infty$,

$$P\{X(t+u\Delta_c) > c - \gamma/c, X(t) < c - y/c\}$$

$$\leq (1+o(1))\psi(c) \int_y^\infty e^z P\{W_t(u) > z - \gamma\} \, dz$$

$$\leq \psi(c) h(y),$$

where $h(y) \to 0$ as $y \to \infty$, establishing (A3). To show that (A5) holds, note that

$$P\{X(t) > c, X(t+u\Delta_c) > c\}$$

$$\leq P\{X(t) + X(t+u\Delta_c) > 2c\}$$

$$\sim \psi\left(\left[\frac{2c^2}{1+\rho(t,t+u\Delta_c)}\right]^{1/2}\right)$$

$$= \psi(c)\left(\frac{1+\rho(t,t+u\Delta_c)}{2}\right)^{1/2} \exp\left[-\frac{c^2}{1+\rho(t,t+u\Delta_c)} + \frac{c^2}{2}\right]$$

$$\leq \psi(c) \exp\left[-\frac{c^2}{2}\left(\frac{1-\rho(t,t+u\Delta_c)}{2}\right)\right].$$

By (2.7), there exists $\eta > 0$ such that $c^2[1 - \rho(t, t+u\Delta_c)] \geq \eta \|u\|^\alpha L(\|u\|)$ for all $t, t+u\Delta_c \in [D]_\delta$. It then follows from (5.37) that (A5) holds with $f(u) = B_\lambda \exp(-u^\lambda)$ with $0 < \lambda < \alpha$, for some $B_\lambda > 0$.

To prove (A4), we use a technique of Fernique [16]. Let $a > 0$, $0 < \zeta < \alpha$, $1 \leq \xi < 2^{\zeta/2}$, $\kappa = \sum_{r=0}^\infty \xi^{-r}$ and $w_r = \xi^{-r}/2\kappa$. Define

(5.38)
$$B_r = \{t + k2^{-r} a\Delta_c : 0 \leq k_i < 2^r, k_i \in \mathbf{Z}\},$$

$$F = \left\{\sup_{u \in I_{t,a\Delta_c}} X(u) > c\right\},$$

$$E_{-1} = \{X(t) \leq c - \gamma/c\},$$

$$E_r = \left\{\sup_{v \in B_r} X(v) \leq c - \gamma(1 - w_0 - \cdots - w_r)/c\right\} \qquad \text{for } r \geq 0,$$

recalling that $\sum_{r=0}^{\infty} w_r = \frac{1}{2}$. Note that $B_r \subset B_{r+1} \subset I_{t,a\Delta_c}$ and that by the continuity of $X$, $P(F \cap E_{-1}) \leq \sum_{r=0}^{\infty} P(E_{r-1} \cap E_r^c)$. Moreover,

$$\begin{aligned}(5.39)\quad &P(E_{r-1} \cap E_r^c) \\ &\leq 2^{r+d} \sup_{\substack{v \in I_{t,a\Delta_c'} \\ \varepsilon \in \{0,1\}^d \setminus \{\mathbf{0}\}}} P\{X(v) \leq c - \gamma(1 - w_0 - \cdots - w_{r-1})/c, \\ &\qquad X(v + \varepsilon 2^{-r} a\Delta_c) > c - \gamma(1 - w_0 - \cdots - w_r)/c\}.\end{aligned}$$

Given $X(v) = c - y/c$, the conditional distribution of $c[X(v + \varepsilon 2^{-r} a\Delta_c) - X(v)]$ is normal with mean $-c(c - y/c)[1 - \rho(v, v + \varepsilon 2^{-r} a\Delta_c)] < 0$ and variance $c^2[1 - \rho^2(v, v + \varepsilon 2^{-r} a\Delta_c)]$, which is bounded by $B(a 2^{-r})^\zeta$ for some $B > 0$, in view of (2.7). Hence

$$\begin{aligned}(5.40)\quad &P\left\{\sup_{\varepsilon \in \{0,1\}^d} c[X(v + \varepsilon 2^{-r} a\Delta_c) - X(v)] > w_r y | X(v) = c - y/c\right\} \\ &\leq 2^d \exp[-C(w_r y)^2/(a 2^{-r})^\zeta]\end{aligned}$$

for some $C > 0$. Similarly, $X(v + \varepsilon 2^{-r}) - X(v)$ has mean 0 and variance bounded by $B'(a 2^{-r})^\zeta/c^2$ for some $B' > 0$ in view of (2.7). Hence by choosing $C$ small enough,

$$\begin{aligned}(5.41)\quad &P\left\{\sup_{\varepsilon \in \{0,1\}^d} c[X(v + \varepsilon 2^{-r} a\Delta_c) - X(v)] > \beta w_r\right\} \\ &\leq 2^d \exp[-C(w_r \beta)^2 c^2/(a 2^{-r})^\zeta].\end{aligned}$$

Let $\eta := 2^\zeta/\xi^2 > 1$. Combining (5.39)–(5.41) with $P\{X(v) \in c - dy/c\} \sim \psi(c) e^y \, dy$ then yields

$$\begin{aligned}(5.42)\quad &P(F \cap E_{-1}) \\ &\leq (1 + o(1))\psi(c) \sum_{r=0}^{\infty} 2^{-r} \int_{\gamma/2}^{\infty} \exp[y - C\eta^r y^2/(4a^\zeta \kappa^2) + C'r] \, dy \\ &\quad + \sum_{r=0}^{\infty} 2^{-r} \exp[-C\eta^r \beta^2 c^2/(4a^\zeta \kappa^2) + C'r]\end{aligned}$$

for some $C' > 0$. Let $\gamma_a = a^{\zeta/3}$ and take $\beta^2 > (2a^\zeta \kappa^2)/C + \lambda$ with $\lambda > 0$. Then for all large $c$ and $\gamma_a \leq \gamma \leq c$, (5.42) is bounded above by $\psi(c) N_a(\gamma)$, where

$$\begin{aligned}N_a(\gamma) = 2 \sum_{r=0}^{\infty} 2^{-r} \bigg\{ &\int_{\gamma/2}^{\infty} \exp[y - C\eta^r y^2/(4a^\zeta \kappa^2) + C'r] \, dy \\ &+ \exp[-C\eta^r \lambda \gamma^2/(4a^\zeta \kappa^2) + C'r] \bigg\}\end{aligned}$$

satisfies $N_a(\gamma_a) + \int_1^\infty y^s N_a(\gamma_a + y)\,dy = o(a^p)$ for all $s > 0$ and $p > 0$. □

PROOF OF THEOREM 4.1. Take $\delta > 0$ and let $D = \{t : 2\delta \leq t_i \leq 1 - 2\delta$ for all $i\}$,

$$\tau(t) = \{F(t)(1 - F(t))\}^{1/2}, \qquad \tau^* = \inf_{t \in [D]_\delta} \tau(t),$$

(5.43)
$$X_c(t) = Z_{n_c}(t)/\tau(t), \qquad b_c(t) = c/\tau(t).$$

Note that $X_c(t)$ has mean 0 and unit variance. We now show that conditions (C) and (B1)–(B5) with $D_c = D$ hold for $X_c(t)$. In view of $F(t + u) = F(t) + u'\nabla F(t) + o(\|u\|)$ and a similar Taylor expansion for $\tau(t + u)$,

$$\rho_c(t, t+u) = F(t)(1 - F(t+u))/\{\tau(t)\tau(t+u)\}$$
$$= 1 - (1 + o(1))u'\nabla F(t)/\{2F(t)(1 - F(t))\}$$

as $u \to 0$. Hence (C) is satisfied with $L(\|u\|) \equiv 1$, $\alpha = 1$, $r_t(u) = u'\nabla F(t)/\{2F(t)(1 - F(t))\}$ and $\Delta_c = (2c^2)^{-1}$. Since $c = o(n_c^{1/6})$ and $\sqrt{n_c} Z_{n_c}(t)$ is a sum of i.i.d. bounded random variables, (B1) holds by moderate deviations theory (cf. [15], Theorem 16.7.1). Moreover, conditioned on $F_{n_c}(t) = x$, $F_{n_c}(t + u\Delta_z) = x + W/n_c$, where $W$ is a Binomial$(n, p)$ random variable with $n = n_c - x$ and $p = \{F(t + u\Delta_z) - F(t)\}/(1 - F(t))$. Making use of this and the functional central limit theorem, it can be shown that (B2) holds and

(5.44)
$$E\{z[X_c(t + u\Delta_z) - X_c(t)] | X_c(t) = z - y/z\}$$
$$\to -u'\nabla F(t)/\{4F(t)(1 - F(t))\},$$

(5.45)
$$\operatorname{Cov}(z[X_c(t + u\Delta_z) - X_c(t)], z[X_c(t + v\Delta_z) - X_c(t)] | X_c(t) = z - y/z)$$
$$\to \sum_{i=1}^d \min(u_i, v_i) \frac{\partial F}{\partial t_i}(t) / \{2F(t)(1 - F(t))\},$$

uniformly over bounded, nonnegative values of $y$ and over $t \in [D]_\delta$ and $c/2 \leq z \leq c/\tau^*$. Note that (5.44) and (5.45), which are analogous to (3.12) and (3.13), give the mean and covariance functions of the limiting Gaussian process $W_t(u)$ in (B2) and are in agreement with the $\alpha$ and $r_t$ of condition (C). The proof of (B3) and (B5) uses ideas similar to those in the proofs of Lemma 3.6(ii) and Lemma 3.7, together with large deviation (instead of Berry–Esseen) bounds for sums of i.i.d. bounded random variables.

To prove (B4), we modify the preceding proof of (A4) in Theorem 2.1 as follows. Let $a > 0$, $1 < \xi < \sqrt{2}$, $\kappa = \sum_{r=0}^\infty \xi^{-r}$, $w_r = \xi^{-r}/2\kappa$ and $c/2 \leq z \leq c/\tau^*$. Pick $r_z$ such that $\theta \leq 2^{r_z} z n_c^{-1/2} < 2\theta$, in which $\theta$ will be specified below. For $u \in I_{v - a2^{-r_z}\Delta_z, a2^{-r_z}\Delta_z}$,

(5.46)
$$Z_{n_c}(v) \geq Z_{n_c}(u) - n_c^{1/2}[F(v) - F(v - a2^{-r_z}\Delta_z)]$$
$$\geq Z_{n_c}(u) - \omega z^{-1},$$

where $\omega > 0$ can be made arbitrarily small by choosing $\theta$ large since $n_c^{1/2} 2^{-r_z}$ lies between $\theta^{-1} z$ and $\theta^{-1} z/2$. Making use of (5.46), we can choose $\theta$ large enough such that

$$
(5.47) \quad \left\{ \sup_{u \in I_{v-a2^{-r_z}\Delta_z, a2^{-r_z}\Delta_z}} z[X_c(u) - X_c(v)] > \gamma/2 \right\} \cap \{X_c(v) \leq z - \gamma/(2z)\} = \varnothing.
$$

For fixed $z$ and $r \geq 0$, define $F, E_{-1}, B_r$ and $E_r$ ($r \geq 0$) by (5.38) in which $c$ is replaced by $z$ and $X(\cdot)$ by $X_c(\cdot)$. Then $P(F \cap E_{-1}) \leq P(F \cap E_{r_z}) + \sum_{r=0}^{r_z} P(E_{r-1} \cap E_r^c)$. We can then proceed as in the preceding proof of Theorem 2.1, using bounds for binomial (instead of normal) tail probabilities.

Verification of (C) and (B1)–(B5) enables us to apply Corollary 2.9 after introducing the change of variables $t \to F(t)$ so that $F$ is a distribution function on $[0,1]^d$ (see the paragraph preceding Theorem 4.1). By (5.44), (5.45) and Lemma 2.3, $H(t) = \{\prod_{i=1}^{d}(\partial F/\partial t_i)(t)\}/\{4F(t)(1-F(t))\}^d$. Moreover, $\inf_{t \in D} 1/\{F(t)(1-F(t))\}^{1/2}(= b_D) = 2$ when $\delta$ is sufficiently small, and $|\nabla_{\perp}^2 b(t)| = \|\nabla F(t)\|^2 \{F(t)(1-F(t))\}^{-3/2}$ and $F(t) = 1/2$ for all $t \in \mathcal{M}$. Hence, applying Corollary 2.9 to $D$ and then letting $\delta \to 0$, we obtain (4.4). Since the probability of joint occurrence of $\{\sup_t Z_{n_c}(t) > c\}$ and $\{\inf_t Z_{n_c}(t) < -c\}$ is negligible compared to (4.4), (4.5) follows from (4.4). □

DEPARTMENT OF STATISTICS
AND APPLIED PROBABILITY
NATIONAL UNIVERSITY OF SINGAPORE
REPUBLIC OF SINGAPORE 119260
E-MAIL: stachp@nus.edu.sg

DEPARTMENT OF STATISTICS
STANFORD UNIVERSITY
STANFORD, CALIFORNIA
USA
E-MAIL: lait@stat.stanford.edu